\documentstyle[11pt,french,amssymb,amsfonts,amsmath,xypic]{article}
\def\tanmoytranslate#1{\catcode`#1=\active\mytanmoytranslate}
\def\mytanmoytranslate#1#2{\def#1{#2}}
\tanmoytranslate    \ä ä{\"a}
\tanmoytranslate    \à à{\`a}
\tanmoytranslate    \á á{\'a}
\tanmoytranslate    \ã ã{\~a}
\tanmoytranslate    \â â{\^a}
\tanmoytranslate    \ë ë{\"e}
\tanmoytranslate    \è è{\`e}
\tanmoytranslate    \é é{\'e}
\tanmoytranslate    \ê ê{\^e}
\tanmoytranslate    \ï ï{\"\i}
\tanmoytranslate    \ì ì{\`\i}
\tanmoytranslate    \í í{\'\i}
\tanmoytranslate    \î î{\^\i}
\tanmoytranslate    \ö ö{\"o}
\tanmoytranslate    \ò ò{\`o}
\tanmoytranslate    \ó ó{\'o}
\tanmoytranslate    \õ õ{\~o}
\tanmoytranslate    \ô ô{\^o}
\tanmoytranslate    \ü ü{\"u}
\tanmoytranslate    \ù ù{\`u}
\tanmoytranslate    \ú ú{\'u}
\tanmoytranslate    \û û{\^u}
\tanmoytranslate    \Ä Ä{\"A}
\tanmoytranslate    \À À{\`A}
\tanmoytranslate    \Á Á{\'A}
\tanmoytranslate    \Ã Ã{\~A}
\tanmoytranslate    \Â Â{\^A}
\tanmoytranslate    \Ë Ë{\"E}
\tanmoytranslate    \È È{\`E}
\tanmoytranslate    \É É{\'E}
\tanmoytranslate    \Ê Ê{\^E}
\tanmoytranslate    \Ï Ï{\"I}
\tanmoytranslate    \Ì Ì{\`I}
\tanmoytranslate    \Í Í{\'I}
\tanmoytranslate    \Î Î{\^I}
\tanmoytranslate    \Ö Ö{\"O}
\tanmoytranslate    \Ò Ò{\`O}
\tanmoytranslate    \Ó Ó{\'O}
\tanmoytranslate    \Õ Õ{\~O}
\tanmoytranslate    \Ô Ô{\^O}
\tanmoytranslate    \Ü Ü{\"U}
\tanmoytranslate    \Ù Ù{\`U}
\tanmoytranslate    \Ú Ú{\'U}
\tanmoytranslate    \Û Û{\^U}
\tanmoytranslate    \ñ ñ{\~n}
\tanmoytranslate    \Ñ Ñ{\~N}
\tanmoytranslate    \ç ç{\c c}
\tanmoytranslate    \Ç Ç{\c C}
\tanmoytranslate    \ß ß{\ss}
\tanmoytranslate    \Æ Æ{\AE}
\tanmoytranslate    \æ æ{\ae}
\tanmoytranslate    \Å Å{\AA}
\tanmoytranslate    \å å{\aa}
\tanmoytranslate    \© ©{\copyright}
\tanmoytranslate    \£ £{\pounds}
\tanmoytranslate    \¶ ¶{\P}
\tanmoytranslate    \§ §{\S}
\tanmoytranslate    \¿ ¿{?`}
\tanmoytranslate    \¡ ¡{!`}

\begin{document}

\author{Gentiana Danila}
\title{Sections du fibré déterminant
sur l'espace de modules des faisceaux semi-stables de rang 2 sur le
plan projectif}
\date{1 Avril 1999}
\maketitle
\newtheorem{theor}{Théorème}[section]
\newtheorem{prop}[theor]{Proposition}
\newtheorem{cor}[theor]{Corollaire}
\newtheorem{lemme}[theor]{Lemme}
\newtheorem{slemme}[theor]{Sous-lemme}
\newtheorem{defi}[theor]{Définition}
\newtheorem{conj}[theor]{Conjecture}
\newtheorem{rem}[theor]{Remarque}
\newtheorem{rems}[theor]{Remarques}
\newtheorem{nota}[theor]{Notations}

\def\fs{{faisceau }}  
\def\fx{{faisceaux }}
\def\alg{{algébrique }}
\def\algs{{algébriques }}
\def\th{{théorème }}
\def\rep{{représentation }}
\def\reps{{représentations }}
\def\irr{{irréductible }}
\def\irrs{{irréductibles }}
\def\coh{{cohomologie }}
\def\co{{cohérent }}
\def\fib{{fibré }}
\def\fibs{{fibrés }}
\def\mor{{morphisme }}
\def\isom{{isomorphisme }}
\def\iso{{isomorphisme }}
\def\mors{{morphismes }}
\def\sur{{surjectif }}
\def\diff{{différentielle }}
\def\diffs{{différentielles }}
\def\inve{{inversible }}
\def\sct{{section }}
\def\hol{{holomorphe }}
\def\hols{{holomorphes }}
\def\app{{application }}
\def\apps{{applications }}

\def\fil{{filtration }}
\def\fils{{filtrations }}
\def\vs{\vskip 4mm}

\def\proj{{\mathbb{P}}}
\def\xm{X^m}
\def\xms{X^m_*}
\def\smx{\es^m(X)}
\def\smxs{\es^m_*(X)}
\def\xim{\Xi_{m}}
\def\pp{\proj_2}

\def\Hilb{{\rm Hilb}}
\def\hil{\Hilb^m(\pp)}
\def\hils{\Hilb^m_*(\pp)}
\def\sm{\es^m(\pp)}
\def\sms{\es^m_*(\pp) }
\def\hilx{\Hilb^m(X)}
\def\hilxs{\Hilb^m_*(X)}
\def\ppm{\pp^m}   
\def\ppms{\proj_{2*}^m}


\def\de{{\mathfrak{d}}}   
\def\vk{\V_k}
\def\vkde{\vk\otimes\de}
\def\ID{{\I_D}}

\def\sigm{{\mathfrak{S}}_m}
\def\comp{{\Bbb C}}
\def\sigmp{{\mathfrak{S}}_{m+1}}

\def\en{{\Bbb N}}
\def\zed{{\Bbb Z}}

\def\rightto#1{\smash{\mathop{\longrightarrow}\limits^{#1}}}
\def\downto#1{\Big\downarrow\rlap{$\vcenter{\hbox{$\scriptstyle#1$}}$}}

\def\ra{\rightarrow}
\def\surto{\twoheadrightarrow}

\def\tens{\otimes}

\def\L{{\mathcal L}}
\def\B{{\mathcal B}}
\def\A{{\mathcal A}}
\def\K{{\mathcal K}}
\def\has{{\mathcal H}}
\def\I{{\mathcal I}}
\def\N{{\mathcal N}}

\def\D{{\mathcal D}}
\def\R{{\mathcal R}}
\def\I{{\mathcal I}}
\def\O{{\mathcal O}}
\def\V{{\mathcal V}}
\def\J{{\mathcal J}}
\def\maU{{\mathcal U}}
\def\maV{{\mathcal V}}
\def\maF{{\mathcal F}}
\def\maG{{\mathcal G}}
\def\h{{\mathfrak{h}}}
\def\g{{\mathfrak{g}}}
\def\uHom{\underline{{\rm Hom}}}

\def\etab{\bar{\eta}}
\def\tL{\widetilde{L}}
\def\tV{\widetilde{V}}
\def\tI{\widetilde{I}}
\def\tW{\widetilde{W}}
\def\tK{\widetilde{K}}
\def\tF{\widetilde{\F}}
\def\tM{\widetilde{M}}
\def\tD{\widetilde{D}}
\def\tmL{\widetilde{\L}}

\def\tna{\widetilde{\nabla}}
\def\tnu{\widetilde{\nu}}

\def\M{{\bf M}}

\def\es{{\rm S}}
\def\H{{\rm H}}
\def\Stab{{\rm Stab\,} }
\def\Ker{{\rm Ker}\, }
\def\coker{{\rm coker}\, }
\def\dim{{\rm dim}\, }
\def\codim{{\rm codim}\, }
\def\supp{{\rm supp}\, } 
\def\Hom{{\rm Hom}\, }
\def\Pic{{\rm Pic}\, }

\def\U{{\rm U}}
\def\sl{{\rm SL}}
\def\sym{{\rm Sym\,}}
\def\spec{{\rm Spec\,}}
\def\F{{\rm F}}      
\def\gr{{\rm gr}}
\def\Im{{\rm Im}}
\def\Tor{{\rm Tor}}
\def\det{{\rm det}}


{\bf Abstract~:} {\small  We provide supporting examples to Le Potier's Strange
duality conjecture, in the case of the moduli space ${\bf M}$ of rank
$2$ semi-stable sheaves on the projective plane, with even first Chern
class, and second Chern class less or equal to $19$. We compute in
this case the dimension of the space of global sections of the
determinant bundle on ${\bf M}$.}



Running heads: Sections du fibré déterminant pour le plan projectif

\section{Introduction}
\label{intro}

Soit ${\bf M}$  l'espace de modules des faisceaux semi-stables 
de
rang $2$, et classes de Chern $c_1=0$, $c_2=n$ sur le plan
projectif complexe. On sait que ${\bf M}$ est une variété projective
irréductible de dimension $4n-3$.

Soit $F$ un \fs semi-stable de
rang $2$, et classes de Chern $c_1=0$, $c_2=n$ sur le plan
projectif. On dit que la droite $H$ est de saut pour $F$ si la 
restriction de $F$ à $H$ n'est pas triviale. On montre que 
(\cite{LeP-Durham},4.1, \cite{Barth}) l'ensemble $\gamma_F$ des
droites de saut pour $F$ est le support d'un \fs cohérent $\Theta$ sur
le plan projectif dual (l'espace des droites de $\pp$). Ce \fs est pur
de dimension $1$ et son support schématique est une courbe  de degré
$n$, appelé courbe des droites de saut de $F$.  
L'application qui associe à la classe de $F$ la courbe $\gamma_F$
définit un \mor $\gamma$ appelé \mor des droites de saut
$F\mapsto\gamma_F$: 
$$\gamma:{\bf  M}\ra  \proj(\H^0(\pp^*,{\cal 
O}_{\proj_{2}^*}(n)))=\proj_N$$
Le fibré déterminant de Donaldson sur ${\bf M}$, noté ${\mathcal D}$, est 
isomorphe à  l'image réciproque $\gamma^*(\O(1))$ de $\O_{\proj_N}(1)$
par $\gamma$ .
On en déduit un  morphisme non nul $\gamma^*$ sur les espaces de sections 
globales :
$$\gamma^*:\H^0(\pp^*,\O(n))^*\ra \H^0({\bf M},{\mathcal D})$$
si on tient compte du fait que 
$\H^0(\proj_N,\O(1))\simeq\H^0(\pp^*,\O(n))^*.$

Le but de cet article est de démontrer le théorème suivant~:
\begin{theor}
\label{teore}
Si $2\le n=c_2\le 19$, l'application linéaire canonique
$$\gamma^*:\H^0(\pp^*,\O(n))^*\ra \H^0({\bf M},{\mathcal D})$$
est un isomorphisme.
\end{theor}

Cet énoncé répond partiellement à une question posée par
A. Beauville, selon laquelle $\gamma^*$ serait un isomorphisme pour
tout $n\in{\mathbb N}$. Il fournit des exemples pour la conjecture de
``Dualité étrange'' de Le Potier.

Remarquons tout d'abord que, comme les deux membres sont des
représentations de ${\sl (3)}$, et que le premier est une 
représentation 
irréductible, le morphisme est injectif, puisque équivariant et
non nul. Il suffit  
donc de calculer la dimension du membre de droite pour conclure. Il suffit en
outre de se placer dans le cas $n\ge 3$, puisque pour $n=0$, l'espace
${\bf M}$  est réduit à un point, pour $n=1$, il est vide, et pour 
$n=2$,
$\gamma$ est un isomorphisme \cite{Barth2}.

La structure de la démonstration est la suivante. On se fixe un entier
positif $l$. On introduit la notion de système cohérent, qui
consiste à considérer en même temps que le \fs $F$, un sous-espace
vectoriel $\Gamma$ de son espace de sections $\H^0(F)$. La dimension
de $\Gamma$ donne l'ordre du système cohérent. À l'aide de
résultats de Min He (\cite{He}) sur les espaces de modules de systèmes cohérents $(\Gamma,F(l))$
d'ordre $1$, dont le \fs sous-jacent est de rang $2$, et de classes de
Chern $c_1=2l, c_2=n+l^2$,
on se ramène au paragraphe \ref{syst}, pour $n$ compris entre $l(l-1)$
et $(l+1)(l+2)$, à l'étude de l'espace des
sections d'un  fibré vectoriel $\es^l\R\tens \de$
sur un ouvert $U$ du schéma de Hilbert $\hil$ des sous-schémas
finis de longueur $m=n+l^2$. Si
$\Xi\subset\hil\times\pp$ est le sous-schéma universel,
 ${\I}$ est le faisceaux d'idéaux associé,
 $pr_1:\hil\times\pp\to \hil$, $pr_2:\hil\times\pp\to \pp$ sont les deux projections, le 
faisceau
algébrique cohérent $\R$ est défini par
$\R=R^1pr_{1*}({\I}(2l-3))$. Ce faisceau  est localement libre 
en de\-hors du fermé de Brill-Noether $B$ des schémas $Z\in 
\hil$
tels que $h^0(I_Z(2l-3))\ne 0$. On note $U$ l'ouvert  complémentaire de
$B$. La codimension de $B$ est supérieure
ou égale à $2$, donc les résultats de cohomologie locale  nous
permettent de passer de $\hil$ à $U$ pour le calcul d'un espace de sections. Le fibré $\de$ est le \fib déterminant sur le schéma de Hilbert
 $\hil$. On désigne par $E$ l'espace de sections $\H^0(\pp,\O(1))$.
Au paragraphe \ref{resolution} on montre que  $\es^l\R\tens \de$ admet
sur $U$ une  résolution par un complexe
 $K^{i}=\Lambda^{-i}\es^kE\tens\es^{l+i}\vkde$ pour $i=0,\ldots,l$,
 où $k=2l-3$,  et
 $\vk$ est défini par $\vk=pr_{1*}(\O_{\Xi}\tens pr_2^*(\O(k)))$. Par
 conséquent, la  suite spectrale
$E^{p,q}_2=\H^q(K^{-p})$ admet pour  aboutissement en degré $0$ l'espace
$\H^0(\es^l\R\tens\de)$. Le tableau suivant représente les termes
$E^{p,q}_2$ de la suite spectrale~:

\begin{tabular}{c|c|c|c|c|c|c}
&&&&&&\\
&$\H^l(\Lambda^l\es^kE\tens\de)$&  & & &&$q=l$\\
&&&&&&\\ \hline
&&&&&&\\
&$\H^{l-1}(\Lambda^l\es^kE\tens\de)$
&$\H^{l-1}(\Lambda^{l-1}\es^kE\tens\V_k\tens\de)$ & & &&$q=l-1$\\
&&&&&&\\ \hline
&&&&&&\\
&$\vdots$&$\vdots$&$\ddots$& &&$\vdots$\\
&&&&&&\\ \hline
&&&&&&\\
&$\H^1(\Lambda^l\es^kE\tens\de)$&$\H^1(\Lambda^{l-1}\es^kE\tens\V_k\tens\de)$&$\cdots$&$\cdots$&$\H^1(\es^{l-1}\V_k\tens\de)$&$q=1$\\
&&&&&&\\
\hline
&&&&&&\\
&$\H^0(\Lambda^l\es^kE\tens\de)$&$\H^0(\Lambda^{l-1}\es^kE\tens\V_k\tens\de)$&$\cdots$&$\H^0(\es^{l-1}\V_k\tens\de)$&$\H^0(\es^l\V_k\tens\de)$&$q=0$\\
&&&&&&\\
\hline
&&&&&&\\
&$p=-l$&$p=-l+1$&$\ldots$&$p=-1$&$p=0$&\\
&&&&&&
\end{tabular}

Nous allons prouver le \th suivant~:
\begin{theor}
\label{enunt} 
On a sur $\hil$~:

i) $\H^0(\de)=E$ et $\H^q(\de)=0$ pour $q>0$;

ii) $\H^0(\vkde)=\es^{2l-2}E\tens\es^{m-1}E$ et
$\H^q(\vkde)=0$ pour $q>0$;

iii) $\dim\H^0(\es^2\V_k\tens\de)=\dim(\es^{2k+1}E\tens\es^{m-1}E\oplus\es^2(\es^{k+1}E)\tens\es^{m-2}E)-\dim(\es^{2k+2}E\tens\es^{m-3}E)$ et $\H^1(\es^2\V_k\tens\de)=0$;

iv) \begin{eqnarray*}
\dim\H^0(\es^3\V_3\tens\de)&=& \dim\es^{10}E\tens\es^{n+8}E+\dim\es^7E\tens\es^4E\tens\es^{n+7}E-\dim\es^{10}E\tens E\tens\es^{n+7}E+\\
&+&\dim\es^{n+6}E(\dim\es^{6,6}E+\dim\es^{7,4,1}E+\dim\es^{8,2,2}E+\dim\es^{6,4,2}E+\dim\es^{4,4,4}E).
\end{eqnarray*}
\end{theor}

Le (i) est une conséquence du \th de Kawamata-Viehweg (\cite{C-K-M}),
expliquée dans le lemme \ref{Kawamata}. Le (ii) fait
l'objet de l'article \cite{D}. On démontre (iii) au paragraphe \ref{calco}, après avoir introduit
 une description appropriée de l'ouvert $\hils$
du schéma de Hilbert $ \hil$, formé par les schémas avec au 
plus un point multiple, qui soit double. Le point le plus délicat  reste
(iv), et il constitue l'objet de la dernière partie. Ce groupe est
obtenu comme noyau d'un \mor de \reps de  ${\sl (3)}$, que nous
expliciterons.

Ces données suffisent pour conclure dans les cas $l=2,3$.
On est ainsi en mesure de calculer la dimension de l'espace $\H^0({\M},{\D})$, dans le cas où $n\le 19$. Pour étendre ce
résultat au cas $n\ge 20$, on a besoin d'étendre le \th d'annulation
de la cohomologie supérieure des \fibs $\es^k\V\tens\de$ sur le schéma
de Hilbert $\hil$ et d'améliorer la méthode utilisée pour le calcul de
$\H^0(\es^3\V_3\tens\de)$, afin de réussir à calculer $\H^0(\es^l\vkde)$, pour $l\ge 3$.

\section{Préliminaires}
\label{prelim}

Notation~: Le corps de base est $\comp$. Pour un espace vectoriel $V$ nous noterons $\proj(V)$ l'espace
projectif des droites de $V$ et $\proj^{.}(V)$ l'espace projectif de
Grothendieck des espaces vectoriels
quotients de dimension $1$.

\subsection{Calculs d'invariants}
\label{calinv}

On considère un ensemble fini $I$ muni d'une action transitive d'un
groupe fini $G$. Soit $Y$ une variété sur laquelle $G$ agit à
gauche. Considérons  pour chaque $i\in I$ un \fib $L_i$ sur $Y$ de fa\c con 
qu'on ait un \iso canonique $h_g:g^*(L_i)\simeq 
L_{g^{-1}(i)}$ pour tout $i\in I$ et $g\in G$, et  pour tous $g,g'\in G$, $h_g\circ h_{g'}=h_{gg'}$. (En 
particulier pour tout $g$ dans $\Stab\{i\}$, le stabilisateur de $i$, on a  $g^*(L_i)\simeq 
L_i$).  On a alors un diagramme commutatif:
$$\begin{array}{ccc}
L_{g^{-1}(i)}&\rightto{h_{g}}&L_i\\
\downto{}&&\downto{}\\
Y&\rightto{g}&Y\end{array}$$

On considère l'espace vectoriel des sections $M_i=\H^0(L_i)$, et la 
somme directe $M=\oplus_{i\in I}M_i$ (espace vectoriel des familles 
$s=(s_i)_i$ avec $s_i\in M_i$).

L'\iso $h_{g}$ induit un \iso $\lambda_{g}:M_i\to 
M_{g(i)}$ en posant pour $x\in Y$
$$\lambda_{g}(s)(x)=h_{g}s(g^{-1}(x)).$$

On peut facilement vérifier que $\lambda_{gg'}=\lambda_{g}\lambda_{g'}$. En 
particulier, ceci fournit une 
action à gauche du stabilisateur de $i$ sur $M_i$. On définit aussi 
une action à gauche de $G$ sur $M$ en posant 
$g(s)_i=\lambda_{g}(s_{g^{-1}(i)})$. Le lemme 
suivant utilisé au paragraphe \ref{calco} et démontré dans \cite{D} est l'ingrédient essentiel 
des calculs d'invariants sur $\xm$:

\begin{lemme}

\label{propinv}
Soit $i\in I$. La projection 
$pr_i:M\to M_i$ induit un \iso $M^G\to M_i^{\Stab\{i\}}$.

\end{lemme}

\subsection{Le gradué d'un produit tensoriel}
\label{gradue}

On considère  deux \fx \algs $\maF$ et $\maG$ sur une variété \alg $Z$,
munis de filtrations
décroissantes
$\maF_{i}$ et $\maG_j$, telles que $\maF_{0}=\maF$ et $\maG_0=\maG$.  On pose  $\gr_{i}\maF=\maF_{i}/\maF_{i+1}.$
Il s'agit de démontrer le lemme suivant:

\begin{lemme}
\label{altalema}
Si 
\begin{eqnarray*}
\underline{\Tor}_1(\gr_i(\maF),\gr_j(\maG))=0
\end{eqnarray*}
 pour tous $i$ et
$j$ alors
pour la filtration associée sur
$\maF \otimes \maG$  on a un isomorphisme canonique
$$ \oplus_{i+j=k} \gr_{i}\maF \otimes \gr_{j}\maG\to \gr_{k}(\maF\otimes \maG).$$
\end{lemme}
 
\par{\bf Preuve~:}

Le problème étant local on peut se placer sur un ouvert affine spectre
 d'un anneau $A$. La donnée de $\maF$ et $\maG$ sur $\spec A$ est
 équivalente à la donnée de deux $A$-modules filtrés $M$ et $N$. On se
 ramène à montrer que si 
\begin{eqnarray*}
\Tor_{1}(\gr_{i}M,\gr_{j}N)=0
\end{eqnarray*}
 pour tous $i$
 et $j$ alors pour la filtration naturelle de
$M \otimes N$  on a un isomorphisme canonique
$$ \oplus_{i+j=k} \gr_{i}M \otimes \gr_{j}N \to \gr_{k}(M \otimes N).$$
On rappelle que la filtration naturelle sur $M\tens N$ est donnée par
$\F^k(M\tens N)=\oplus_{i+j=k}\Im(M_i\tens N_j)$. Ici l'hypothèse
permet d'écrire $\F^k(M\tens N)=\oplus_{i+j=k}M_i\tens N_j$.
 On a évidemment un morphisme canonique
$$ \oplus_{i+j=k} \gr_{i}M \otimes \gr_{j}N\to\gr_{k}(M \otimes N)$$
dont il s'agit de vérifier que c'est un isomorphisme.
On raisonne par récurrence sur $k.$ C'est vrai pour $k=0$ par la proposition 6, page 8, chap. III-1 de \cite{Bourbaki}.
On a  clairement pour $\ell>k$
$$\gr_{k}(M \otimes N)= \gr_{k}(M/M_{\ell} \otimes N)$$
de sorte que l'on peut supposer en remplaçant au besoin $M$ par
$M/M_{\ell}$
 que la filtration de $M$  est finie. On écrit la suite
exacte
$$0\to N_{1}\to N\to \gr_{0}N\to 0$$
et si on tensorise par $M$ il résulte de l'hypothèse  que
$\Tor_{1}(M,\gr_{0}N)=0$ et on a alors  une suite exacte
$$0\to M\otimes N_{1}\to M \otimes N\to M\otimes
\gr_{0}N\to 0.$$   De plus, la filtration induite sur $M\otimes N_{1}$  est
bien (au décalage d'indice près) la filtration obtenue à partir du
produit
tensoriel.  On a alors  un diagramme commutatif  de suites exactes
$$\diagram
0\rto& \oplus_{i+j=k,j\geq 1}\gr_{i}M \otimes \gr_{j}N
\rto\dto&
\oplus_{i+j=k} \gr_{i}M\otimes \gr_{j}N\rto\dto&\gr_{k}M\otimes \gr_{0}N \rto\dto&0\\
0 \rto&\gr_{k-1}(M\otimes N_{1})\rto&\gr_k(M\otimes
N)\rto&\gr_{k}(M\otimes
\gr_{0}(N))\rto& 0 
\enddiagram$$
L'hypothèse de récurrence appliquée à $M$  et $N_{1}$ (qui
satisfont bien
entendu aux mêmes hypothèses que $M$ et $N$) montre que la première
flèche
verticale est un isomorphisme; la dernière l'est par l'hypothèse
$\Tor_{1}(\gr_{i}(M),\gr_{0}(N))=0$. Donc la flèche verticale du milieu est
un
isomorphisme.$\Box$

On peut vérifier que dans les conditions dans lesquelles on va
l'appliquer, les hypothèses de ce lemme sont satisfaites~:

\begin{lemme}
\label{incaolema}
Soient $X$ et $Y$ deux variétés \algs munies des \fx \algs $\maF$ et
$\maG$. Alors
$$\underline{\Tor}_1^{\O_X\boxtimes\O_Y}(\maF\boxtimes\O_Y,\O_X\boxtimes\maG)=0.$$
\end{lemme}

\par{\bf Preuve~:}

Le problème étant local on peut travailler sur des ouverts de $X\times
Y$ qui sont de la forme $\spec A\times\spec B$ pour $A$ et $B$ deux
algèbres de type fini sur le corps de base $K$. La donnée de $\maF$ et
$\maG$ sur $\spec A$, respectivement $\spec B$, est équivalente à la
donnée d'un $A$-module $M$ et un $B$-module $N$. On se ramène à
prouver l'annulation de 
$$\Tor^{A\tens_K B}_1(M\tens_K B,A\tens_K N)=0.$$
À cette fin on considère  $Q^{\cdot}$ une résolution projective de $M$. Alors
$Q^{\cdot}\tens_K B$ est une résolution projective pour $M\tens_K
B$. On tensorise cette résolution par $A\tens_K N$. Le complexe
obtenu est en fait $Q^{\cdot}\tens_K N$ qui est exact.$\Box$

\subsection{Représentations irréductibles de $\sl(E)$}
\label{repirr}

Fixons d'abord quelques notations. Soit $E$ un espace vectoriel de 
dimension $r$. Un diagramme de Young est une configuration de boîtes associée à une suite décroissante 
$\lambda_1\ge\ldots\ge\lambda_l >0$ d'entiers, avec $\lambda_i$ boîtes
dans la $i$-ème colonne, les colonnes étant alignées de gauche  
à droite. Un tableau de Young est un diagramme de Young muni d'une numérotation des cases $(i,j)_{1\le i\le l, 
1\le j\le\lambda_i}$.

Soit $d=|\lambda|:= \sum_i\lambda_i$ le nombre total de boîtes. À un tel tableau de Young est 
associée une représentation irréductible de ${G=\sl(E)}$ de la 
manière suivante: on considère l'action naturelle  de ${\mathfrak
  S}_d$ sur $E^{\tens d}$ par permutation des facteurs, le 
sous-groupe $P_{\lambda}$ qui laisse invariantes les colonnes du tableau 
de Young, et le sous-groupe $Q_{\lambda}$ qui laisse invariantes les 
lignes.  Alors, si on pose $a_{\lambda}=\sum_{g\in P_{\lambda}}g$ et 
$b_{\lambda}=\sum_{g\in Q_{\lambda}}sgn(g)g$ (les sommes sont prises
dans l'algèbre du groupe $G$), la représentation 
$\es^{\lambda}E$ est définie par
$$\es^{\lambda}(E)=b_{\lambda}a_{\lambda}(E^{\tens d}).$$

Le \th 6.3 [F-H], page 77, nous donne la dimension de cette 
représentation:
$$\dim\es^{\lambda}E=\prod_{1\le i<j\le 
r}\frac{\lambda_i-\lambda_j+j-i}{j-i}$$
si $\lambda_1\ge\ldots\ge\lambda_r\ge 0$ et $0$ si $\lambda_{r+1}\ne 0$.

Une autre interprétation, via le \th de Bott, est obtenue en considérant la 
variété $D(E)=Drap(E)$ des drapeaux complets 
$0\subset\F_r\subset\F_{r-1}\subset\ldots\subset\F_1=E$, avec 
$\dim\F_i=r-i+1$, puis sur  $D(E)$ 
les \fibs $Q_i$ de rang $1$ définis par les quotients canoniques $Q_i=\F_i/\F_{i+1}$. 
Alors si $\lambda$ est un tableau de Young, la représentation
$$\H^0(D(E),Q_1^{\lambda_1}\tens\ldots\tens Q_r^{\lambda_r})$$
est une représentation irréductible (et les autres groupes de \coh 
s'annulent, par \cite{Demailly}, \cite{Bott}). C'est la 
représentation $\es^{\lambda}(E)$.

On se servira pour le calcul de $\H^0(\es^3\V_k\tens\de)$ des résultats
suivants, qui constituent un court  
résumé de l'étude présentée dans \cite{F-H}.

\vs Étant donnée une représentation $V$ de ${\sl(3)}$ (ce qui revient au 
même que la donnée d'une représentation  de son algèbre de Lie 
${\mathfrak{sl}}(3)$, puisque ${\sl(3)}$ est un groupe de Lie connexe et 
simplement connexe: \cite{F-H},page 109),  on connaît sa décomposition en 
sous-représentations \irrs. On étudie  pour cela l'ensemble de ses 
poids, qui sont les valeurs propres  pour l'action du sous-espace 
$\h\subset {\mathfrak{sl}}(3)$ des matrices diagonales sur $V$. Ces 
valeurs propres sont des formes linéaires sur $\h$. Explicitement on a
$$\h=\{
\left(
\begin{array}{ccc}
a_1&0&0\\
0&a_2&0\\ 
            0&0&a_3\\
            \end{array}
        \right)
     :a_1+a_2+a_3=0
     \}
$$
et on peut donc écrire $\h^*=\comp[L_1]\oplus\comp[L_2]\oplus\comp[L_3]/(L_1+L_2+L_3)$.

\vs Les poids $L_i-L_j\in\h^*$ sont spéciaux, car ils sont les poids de la 
\rep adjointe, et ils sont appelés racines. Les racines engendrent un 
réseau à l'intérieur de $\h^*$, noté $\Lambda_R$. Les poids 
${\alpha}$ d'une \rep finie quelconque se trouvent sur le réseau 
$\Lambda_W\subset\h^*$ engendré par les $L_i$ et sont congrus modulo 
$\Lambda_R$. En choisissant un ordre sur les $L_i$, par exemple 
$L_1>L_2>L_3$, on peut séparer les $6$ racines en $3$ racines 
positives: ${L_1-L_2,L_1-L_3,L_2-L_3}$ et $3$ négatives. Les espaces 
propres correspondants aux racines positives, $\g_{L_1-L_2}$, 
$\g_{L_1-L_3}$, et $\g_{L_2-L_3}$, sont engendrés par les matrices 
$E_{i,j}=(e_{kl})$, avec un seul terme non nul $e_{ij}=1$ au-dessus 
de la diagonale. Une représentation $\h$ de ${\sl(3)}$ étant donnée, il existe un vecteur propre pour $\h$, $v\in V$, qui 
est annulé par $E_{1,2}$, $E_{1,3}$, et $E_{2,3}$, appelé vecteur de 
plus haut poids. Le poids correspondant est appelé poids dominant. Pour une \rep \irr le vecteur de plus haut poids est unique à multiplication par 
des scalaires près.

\vs L'ensemble des poids est préservé par l'action du groupe des 
symétries par rapport aux droites engendrés par les $L_i$ (le 
groupe de Weyl). On en déduit qu'un poids dominant doit se 
trouver dans le $(\frac{1}{6})$-plan délimité par $L_1$ et 
$L_1+L_2=-L_3$, donc il doit être de la forme $aL_1+bL_2+cL_3$ avec 
$a\ge b\ge c$.

\vs La prop. 12.11 [F-H] nous assure que la sous-\rep $W$ de $V$ engendrée 
par les images de $v= aL_1+bL_2+cL_3$ par  application successive des 
trois opérateurs $E_{1,2}$, $E_{1,3}$, et $E_{2,3}$ est irréductible. C'est la 
\rep $\es^{a,b,c}E$, où $E$ est la \rep standard. On la note aussi 
$\es^{\lambda}E$ où $\lambda$ est la partition $(a,b,c)$. Comme $c$ est 
superflu puisque $ L_1+L_2+L_3=0$, on préfère l'écriture 
réduite $\es^{\lambda}E$ avec $\lambda=(a-c,b-c)$.

\vs Maintenant, la formule de Pieri permet de calculer le produit tensoriel 
d'une \rep irréductible et d'une puissance symétrique. On a en effet
$$\es^{\lambda}E\tens\es^pE=\sum_{\stackrel{\lambda_i\le\mu_i\le\lambda_{i-1
}}{|\mu|=|\lambda|+p}}\es^{\mu}E$$

Par exemple 
$\es^8E\tens\es^4E=\es^{12}E+\es^{11,1}E+\es^{10,2}E+\es^{9,3}E+\es^{8,4}E $.

Par une analyse directe on peut trouver la décomposition de n'importe 
quelle \rep. Par exemple, d'après le programme Lie \cite{LiE}, la décomposition du pléthysme $\es^3(\es^4E)$ est donnée par:

$\es^3(\es^4E)=\es^{12}E+\es^{10,2}E+\es^{9,3}E+\es^{8,4}E+\es^{6,6}E+\es^{7,4,1}
+\es^{8,2,2}E+\es^{6,4,2}E+\es^{4,4,4}E$.

\goodbreak

\subsection{Contractions}

\begin{prop}

\label{contr}

Soit $F$ un espace vectoriel. Soient $a$ et $b$ deux entiers positifs, 
tels que $a\ge b$. Alors le \mor canonique (appelé par [F-H] \mor de 
contraction, p.83)
$$c:\es^aF\tens\es^b F\to\es^{a+1}F\tens\es^{b-1}F$$
défini par $y\tens x_1\cdots x_b\mapsto\sum_i 
yx_i\tens x_1\cdots\check{x_i}\cdots x_b$
est surjectif, de noyau $\es^{a,b}F$.
\end{prop}

\par{\bf Preuve:} À cause de la décomposition de Pieri en somme directe de facteurs irréductibles, il suffit de montrer la surjectivité. On se place sur l'espace 
projectif de Grothendieck $\proj^{\cdot}(F)$ et on écrit la suite 
exacte d'Euler associée
$$0\to U\to F\tens\O\to\O(1)\to 0$$
où $U$ est le sous-\fib canonique du \fib trivial de fibre $F$, et 
$\O(1)$ le \fib inversible quotient canonique. On en déduit une suite exacte 
pour les puissances symétriques
$$0\to\es^bU\tens\O(a)\to\es^bF\tens\O(a)\stackrel{\alpha}{\to}\es^{b-1}F\tens\O(a+
1)\to 0$$

On sait que $F=\H^0(\proj^{\cdot}(F),\O(1))$. Le \mor associé à 
$\alpha$ sur les sections globales $\H^0(\alpha)$ coïncide avec $c$, 
donc la surjectivité de $c$ revient à montrer l'annulation de 
l'espace $\H^1(\es^bU(a))$.

\begin{lemme}

On a  $\H^q(\proj^{\cdot}(F),\es^bU(a))=0$ si $a\ge b$, $q\ge 1$ et 
$\dim\, F=3$.

\end{lemme}

\par{\bf Preuve du lemme:}

On considère la variété $D=D(F)$ des drapeaux $h\subset k\subset 
F$ avec $\dim\, h=1$, et $\dim\, k=2$. Soit $\maV\subseteq\maU\subseteq F\tens\O_D$ la \fil 
canonique correspondante 
du \fib trivial de fibre $F$ sur $D$, avec $\maV$ \fib de rang $1$ et 
$\maU$ le sous-\fib canonique de rang $2$. La \fil de $\es^b\maU$ associée 
à la suite exacte
$$0\to\maV\to\maU\to\maU/{\maV}\to 0$$
est donnée par $\F^1\es^b\maU=\es^{b-1}\maU\tens\maV\subset\es^b\maU$.

Si on considère la projection canonique $D\stackrel{pr}{\to}\proj(F)$ 
définie par $(h,k)\mapsto k$, sur la fibre $D_k=\proj(k^*)$ de cette 
projection au-dessus d'un point $k$, le \fib $\maU$ est trivial et le 
\fib $\maV$ s'identifie à $\O(-1)$ de la fibre. On obtient en particulier que $R^qpr_*(\F^1\es^b(\maU))=0$ pour tout 
$q\ge 0$, d'où
$$R^qpr_*(\F^1\es^b(\maU)\tens pr^*(\O(a)))= 
R^qpr_*(\F^1\es^b(\maU))\tens \O(a)=0,$$ 
pour tous $ q\ge 0$ et donc par la 
suite spectrale de Leray on déduit que $\H^q(D,\F^1\es^b(\maU)(a))=0$ 
pour tout $q\ge 0$.

Mais $pr$ est un \mor projectif dont les fibres sont des espaces 
projectifs (c'est un \fib en projectifs) donc $R^q(f_*\O_D)=\O$ si 
$q=0$ et $0$ sinon.

Par suite, comme $\maU=pr^*(U)$ et $\es^b\maU=pr^*(\es^bU)$, on 
obtient 
$$\H^q(\proj(F),\es^bU(a))=\H^q(D,\es^b\maU(a))=\H^q(D,\es^b(\maU/{\maV})(a)),$$
pour tous $ q\ge 0.$ Mais $\O_D(1)=Q_1$ et $\maU/{\maV}=Q_2$ 
sont les quotients canoniques de rang $1$ associés à la \fil 
$\maV\subseteq\maU\subseteq F\tens\O_D$, et cf. [Demailly], on a l'annulation de $\H^q(D,Q_1^{\tens a}\tens Q_2^{\tens b}\tens Q_3^{\tens c})$ pour $q\ge 
1$ et $a\ge b$. D'où le résultat.

\subsection{Application: le noyau et le conoyau du morphisme $\nu$}

Le morphisme $\nu$ est le \mor canonique
$$\nu:\es^3(\es^4E)\to\es^8E\tens\es^4E$$
induit par l'application linéaire $stu\mapsto st\tens u+su\tens
t+ut\tens s$, où $s, t, u \in \es^4E$.

\begin{lemme}
\label{nu}
Le \mor $\nu$ n'est ni injectif, ni surjectif. Son conoyau est 
$\es^{11,1}E$ et son noyau est isomorphe à
$$\es^{6,6}E+\es^{7,4,1}E+\es^{8,2,2}E+\es^{6,4,2}E+\es^{4,4,4}E$$
\end{lemme}
\par{\bf Preuve du lemme:}
On regarde les poids de $\es^4E$. Ils sont donnés par des sommes de $4$ 
poids de $E$, soit: $4L_i$, $3L_i+L_j$, $2L_i+2L_j$, $2L_i+L_j+L_k$, 
avec $i$, $j$, $k$ deux à deux distincts et variant de $1$ à $3$. Les 
poids de $\es^3(\es^4E)$ sont donnés par des sommes de $3$ poids 
$\alpha,\beta,\gamma$ de $\es^4E$. 
On écrit $\alpha+\beta+\gamma$ sous la forme 
$mL_1+nL_2+pL_3$ avec $m+n+p=12$. Si on veut obtenir 
$\alpha+\beta+\gamma=8L_1+4L_2$ on a $m-p=8$, $n-p=4$ donc 
$m+n-2p=m+n+p=12$ donc $p=0$. On ne se sert donc pas des poids qui 
font apparaître $L_3$. Il y a exactement $4$ fa\c cons d'obtenir 
$8L_1+4L_2$:

$4L_1+4L_1+4L_2$ qui correspond au vecteur propre 
$\omega_1=(e_1^4)^2(e_2^4)$,

$4L_1+3L_1+L_2+3L_2+L_1$ qui correspond au vecteur propre 
$\omega_2=(e_1^4)(e_1^3e_2)(e_1e_2^3)$,

$4L_1+2L_1+2L_2+2L_1+2L_2$ qui correspond au vecteur propre 
$\omega_3=(e_1^4)(e_1^2e_2^2)^2$, et

$3L_1+L_2+3L_1+L_2+2L_1+2L_2$ qui correspond au vecteur propre 
$\omega_4=(e_1^3e_2)^2(e_1^2e_2^2)$.

Ces $4$ vecteurs sont indépendants dans $\es^3(\es^4E)$ (ce sont même 
des vecteurs d'une base de $\es^3(\es^4E)$).
On vérifie que leurs images par $\nu$ sont non nulles. Comme $8L_1+4L_2$ figure parmi les poids 
des \reps $\es^{12}E$, $\es^{10,2}E$, $\es^{9,3}E$ et $\es^{8,4}E$, et qu'un 
vecteur de poids $\alpha$ s'envoie sur un vecteur de même poids 
$\alpha$ par le \mor équivariant $\nu$, il résulte que la 
multiplicité de ce poids est $1$ dans chacune de ces représentations, et que 
l'image de $\nu$ atteint chacune de ces \reps vues comme sous-\reps 
\irrs de $\es^8E\tens\es^4E$. $\Box$

\begin{cor}
\label{tnu}
Le \mor $\tnu:\es^3(\es^4E)\tens\es^{m-3}E\to\es^8E\tens\es^4E\tens\es^{m-3}E$ défini par $\tnu=\nu \tens id$ n'est ni injectif, ni surjectif. Son conoyau est 
$\es^{11,1}E\tens\es^{m-3}E$ et son noyau est isomorphe à
$$(\es^{6,6}E+\es^{7,4,1}+\es^{8,2,2}E+\es^{6,4,2}E+\es^{4,4,4}E)\tens\es^{m-3}E$$
\end{cor}
\par{\bf Preuve du corollaire~:}
Comme $\tnu=\nu\tens id$, par le lemme \ref{nu} on obtient le 
résultat.$\Box$



\section{Systèmes cohérents}
\label{syst}

Le but de cette section est de montrer comment on peut ramener le
calcul du nombre de sections du fibré déterminant de Donaldson à 
un
calcul de sections d'un faisceau localement libre sur un ouvert du
schéma de Hilbert. La méthode repose sur un résultat de Min He, 
qui
l'utilisait pour calculer les nombres de Donaldson \cite{He}. On commence par quelques généralités sur les systèmes cohérents.

\subsection{Systèmes cohérents $a$--semi-stables}

On désigne par $K(\pp)$ 
l'algèbre de Grothendieck  des classes de faisceaux algébriques
cohérents sur $\pp$, ou ce qui revient au même,  des classes de fibrés
vectoriels algébriques sur $\pp$.  Cette algèbre est engendrée par la
classe $\eta$ du \fs structural d'une droite. En tant que groupe
abélien, elle est isomorphe à $\zed^3$, un \iso étant donné par le
rang $r$, la classe de Chern $c_1$ et la caractéristique
d'Euler-Poincaré $\chi$. Elle est munie de la forme quadratique
entière non dégénérée définie par
$$2r\chi+c_1^2-r^2.$$ 
La notion d'orthogonalité utilisée par la suite est relative à cette
forme quadratique.

Une classe de Grothendieck $a\in K(\pp)\tens{\Bbb Q}$ est dite {\it
  positive} si le polynôme de Hilbert de $a$ est positif. Étant 
donné
  un faisceau algébrique cohérent $F$ sur $\pp$, on désigne par 
$c(F)$
  la classe de Grothendieck de $F$ dans $K(\pp)$.

\begin{defi}
{\rm Un système cohérent sur $\pp$ est un couple $\Lambda=(\Gamma,F)$ 
où
$F$ est un faisceau cohérent, et $\Gamma$ un sous-espace vectoriel de
$\H^0(F)$. L'ordre du système cohérent est la dimension de l'espace
vectoriel $\Gamma$.}
\end{defi}

Soit $a\in K(\pp)\tens{\Bbb Q}$ une classe de Grothendieck strictement
positive. À un système $\Lambda=(\Gamma,F)$ on associe la classe 
de
Grothendieck $c_a(\Lambda)\in K(\pp)\tens{\Bbb Q}$ définie par
$$c_a(\Lambda)=\dim\Gamma\cdot a+c(F).$$

La catégorie des systèmes cohérents n'est pas une catégorie 
abélienne,
mais elle se plonge dans une catégorie abélienne (la catégorie des
systèmes algébriques) qui a suffisamment d'injectifs. On ne 
considère
ici que des systèmes cohérents $\Lambda=(\Gamma,F)$ dont le 
faisceau
$F$ sous-jacent est de rang $r>0$. Ce rang sera aussi appelé le rang
du système cohérent.

\begin{defi}
{\rm On dit qu'un système cohérent $\Lambda=(\Gamma,F)$ est 
$a$--semi-stable
si
\begin{itemize}
\item{\rm (i)} le faisceau $F$ est sans torsion;
\item{\rm (ii)} pour tout sous-faisceau cohérent $F'\subset F$ de rang
  $r'>0$ on a dans $K(\pp)\tens{\Bbb Q}$:
$$\frac{c_a(\Lambda')}{r'}\le\frac{c_a(\Lambda)}{r}.$$
où $\Lambda'$ est le système cohérent $\Lambda'=(\Gamma',F')$ 
défini
par $\Gamma'=\Gamma\cap \H^0(F')$.
\end{itemize}}
\end{defi}

Dans le cas particulier où $c$ est de rang $2$, et $\dim \Gamma=1$, 
seul cas
utile dans la suite, l'inégalité (ii) signifie que pour tout
sous-faisceau cohérent $F'\subset F$ de rang $1$ on a
\[c(F')\le\left\{ \begin{array}{ll}
        \frac{1}{2}(c(F)-a) & \mbox{si $\Gamma\subset \H^0(F')$}\\
                            &\\
        \frac{1}{2}(c(F)+a) & \mbox {sinon}
                 \end{array}
 \right. \]

\subsection{L'espace de modules $Syst_a(c,k)$}
Soit $c\in K(\pp)$ une classe de Grothendieck fixée et $k$ un entier
$\ge 0$; il existe un espace de modules grossier de systèmes cohérents
$a$--semi-stables $\Lambda=(\Gamma,F)$ tels que $c(F)=c$, et $\dim
\Gamma =k$: c'est une variété projective qui sera notée
$Syst_a(c,k)$. Lorsque $a$ varie, la structure de l'espace de modules
grossier 
$Syst_a(c,k)$ change au passage de certaines valeurs de $a$ appelées
valeurs critiques.

\subsection{Valeurs critiques}

Soit $F$ un faisceau algébrique cohérent de rang $2$, de classes de
Chern $c_1$ et $c_2$ sur le plan projectif. On désigne par $\eta$ la
classe dans $K(\pp)$ du faisceau structural d'une droite. Ainsi, les
faisceaux considérés ont pour classe de Grothendieck
$$c(F)=2+c_1\eta+(\frac{c_1(c_1+1)}{2}-c_2)\eta^2.$$

Soit $l>0$ un entier, fixé dans toute la suite. Les faisceaux
cohérents $F$ de rang $r=2$ et de classes de Chern $c_1=2l$, et
$c_2=n+l^2$ ont alors pour classe de Grothendieck
$$c(l)=2+2l\eta+(l(l+1)-n)\eta^2.$$
On considère l'espace de modules $S_a=Syst_a(c(l),1)$ des classes de
$S$-é\-qui\-va\-lence de systèmes cohérents $a$--semi-stables
$\Lambda=(\Gamma,F)$, d'ordre $1$, où $F$ est un faisceau cohérent 
de
classe de Grothendieck \hbox{$c(F)=c(l)$}. Si cet espace de mo\-dules est
non-vide, une section de $F$ donne une inclusion du faisceau trivial
dans $F$, et donc l'inégalité
$$0<a\le c(l)-2=2l\eta+(l(l+1)-n)\eta^2.$$

\begin{defi}
{\rm Les valeurs critiques pour la famille d'espaces de modules de systèmes
cohérents $S_a$ sont les classes $a\in K(\pp)\tens{\Bbb Q}$ pour 
lesquelles
il existe des systèmes cohérents strictement semi-stables 
relativement
à $a$.}
\end{defi}

Les valeurs critiques sont en nombre fini. On peut en fait les
calculer explicitement, mais on n'aura pas besoin de cela dans la
suite.

\subsection{Les résultats de Min He}

Étant donnée une valeur critique $a$, on désigne par $a_{-}$ et 
$a_{+}$ des
classes de Grothendieck $>0$ encadrant $a$ et telles que dans
l'intervalle $]a_{-},a_{+}[$, $a$ soit la seule valeur critique. On
désigne par $a_{max}=c(l)-2$ la plus grande valeur critique. Pour
$a>a_{max}$, l'espace de modules $S_a$ est vide.

\begin{theor}
\label{minhe}
(Min He) On suppose $n\ge l(l-1)$ et $n\ge 3$.
\begin{itemize}
\item{\rm (i)} L'espace de modules $S_a$ est une variété 
irréductible, normale, de
dimension $\delta=3n+l(l+3)-2$.
\item{\rm (ii)} Si $a$ est valeur critique distincte de $a_{max}$, on 
dispose de
morphismes surjectifs
$$\diagram
S_{a_{-}}\drto_{\pi_{-}}&- - \to&S_{a_{+}}\dlto^{\pi_{+}}\\
&S_a&\\
\enddiagram$$

Au-dessus de l'ouvert des points stables de $S_a$ ces morphismes sont
des isomorphismes. L'image réciproque du fermé $\Sigma$ des points
strictement semi-stables $\Sigma_{-}=\pi^{-1}_{-}(\Sigma)$
(resp. $\Sigma_{+}=\pi^{-1}_{+}(\Sigma)$) est le fermé des points
$a_{+}$-instables (resp. $a_{-}$-instables). Les fermés $\Sigma,
\Sigma_{-}, \Sigma_{+}$ sont de codimension $\ge 2$.
\item{\rm (iii)} Si $a$ n'est pas valeur critique, il existe un 
système
  cohérent universel $\Lambda=({\mathcal V},{\mathcal F})$ 
paramétré par $S_a$.
\end{itemize}
\end{theor}

Ainsi, $S_{a_{+}}$ s'obtient à partir de $S_{a_{-}}$ en rempla\c cant
le fermé $\Sigma_{-}$ par le fermé $\Sigma_{+}$ des points
\hbox{$a_{-}$-instables}.

\subsection{L'espace de modules $S_{a_{max_{-}}}$}

On suppose ici que  $n\ge l(l-1)$ et $n\ge 3$.

L'espace de modules $S_{a_{max}}$ s'identifie au schéma de Hilbert des sous-schémas finis de longueur $n+l^2$ de $\pp$, 
$\Hilb^{n+l^2}(\pp)$. En effet, la condition de semi-stabilité nous
assure que le conoyau du morphisme d'évaluation $\Gamma\tens\O\ra F$
est sans torsion. On connaît la description de tels faisceaux sur
$\pp$. Ils s'écrivent comme $I_Z(c_1)$ où $I_Z$ est l'idéal d'un
sous-schéma fini $Z$ de longueur $c_2$ de $\pp$. De plus l'extension:
$$0\to \Gamma\tens\O\ra F\to I_Z(2l)\to 0$$
donne une filtration de Jordan-Hölder pour $F$, puisque
$a=a_{max}$. Le schéma $Z$ est dans ce cas de longueur
$n+l^2$. L'application qui associe à $F$ le sous-schéma $Z$
correspondant donne l'identification désirée.

Dans le cas où $a=a_{max_{-}}$ on dispose encore d'une extension
$$0\to \Gamma\tens\O\ra F\to I_Z(2l)\to 0$$
soit d'une extension non-triviale de systèmes cohérents
$$0\to\Lambda'\to \Lambda\to\Lambda''\to 0$$
où $\Lambda''=(0,I_Z(2l))$, et $\Lambda'=(\comp,\O_{\pp})$.
Réciproquement, pour $a$ voisin de $a_{max}$, une telle extension
non-triviale définit un système cohérent $a$--semi-stable. Pour
déterminer $S_{a_{max_{-}}}$ il s'agit donc de paramétrer les
extensions non-triviales de ce type.

Mais pour un sous-schéma $Z$, ces extensions sont paramétrées par
$$\proj(Ext^1(I_Z(2l),\O))=\proj^{.}(\H^1(I_Z(2l-3))).$$
La variété
$S_{a_{max_{-}}}$ va donc s'identifier à un fibré en espaces
projectifs associé à un faisceau algébrique cohérent sur
$\Hilb^{n+l^2}(\pp)$ qui ait pour fibre en $Z$, $\H^1(I_Z(2l-3))$. On
obtient le
\begin{theor}
Soit $\Xi\subset\Hilb^{n+l^2}(\pp)\times\pp$ le sous-schéma universel,
et ${\mathcal I}$ le faisceaux d'idéaux associé. Considérons le 
faisceau
algébrique cohérent ${\mathcal R}=R^1pr_{1*}({\mathcal I}(2l-3))$. 
Alors
$$S_{a_{max_{-}}}=\proj({\mathcal R}).$$
\end{theor} 
Ce faisceau ${\mathcal R}$ est localement libre de rang
$$\chi(\O_Z)-\chi(\O(2l-3))=n+1-(l-1)(l-2)$$
en de\-hors du fermé de Brill-Noether $B$ des $Z\in 
\Hilb^{n+l^2}(\pp)$
tels que $h^0(I_Z(2l-3))\ne 0$. On sait d'après Ellingsrud et Str\o mme
(dans \cite{He}, lemme 4.9, en utilisant \cite{E-Stro}, th. 1.1 et cor. 1.2), que si $C$ est une courbe de $\pp$, même non-réduite, le sous-schéma de $\hil$ des points $Z\subset C$ est de dimension $m$. Ceci permet de majorer la dimension du fermé de Brill-Noether $B$ ou de minorer sa codimension~: elle est supérieure à $n-l(l-3)$, donc 
supérieure
à $2l$ si $n\ge l(l-1)$. On peut minorer aussi la codimension de l'image réciproque de 
$B$
dans $\proj({\mathcal R})$~:
\begin{lemme}
Supposons $n\ge l(l-3)$. L'image réciproque de 
$B$
dans $\proj({\mathcal R})$
 est de codimension $\ge n-l(l-1)+2$.
\end{lemme}

\par{\bf Preuve du lemme~:}
Soit $s$ un entier tel que $0<s\le 2l-3$ et $B_s\subset B$ le sous-ensemble localement fermé correspondant aux idéaux $I$ tels que $h^0(I(s-1))=0$, et $ h^0(I(s))\ne 0$. Pour $I\in B_s$ on a 
$$h^1(I(2l-3))\le h^1(I(s-1))=n+l^2-\frac{1}{2}s(s+1).$$
On obtient ainsi une majoration de la dimension de l'image réciproque de $B_s$ par
$$n+l^2+\frac{1}{2}(s+1)(s+2)-1+n+l^2-\frac{1}{2}s(s+1)-1$$
c'est-à-dire
$$2(n+l^2)+s-1.$$
L'assertion résulte en utilisant la majoration $s\le 2l-3$.$\Box$

On note $\pi$ le morphisme canonique de $S_{a_{max_{-}}}$ dans
$\Hilb^{n+l^2}(\pp)$.

\subsection{Le morphisme $f:S_{\epsilon}\to\M_c$}

On désigne par $\M_c$ l'espace de modules des faisceaux
semi-stables de classe de Grothendieck $c$. Ici $c$ est de rang $2$, $c_1=0$ et
$c_2=n$ donc $\M_c=\M$.
Considérons une classe de Grothendieck $\epsilon>0$ inférieure à 
la
plus petite valeur critique. Si $\Lambda=(\Gamma,F)$ est un système
cohérent $\epsilon$--stable, le faisceau $F$ sous-jacent est
semi-stable. Par suite, on obtient un morphisme
$f:S_{\epsilon}\to\M_c$ qui associe à la classe du système 
cohérent
$(\Gamma,F)$ la classe du faisceau $F(-l)$.

\begin{theor}
\label{epsi}
Si $3\le n<(l+1)(l+2)$, le morphisme $f:S_{\epsilon}\to\M_c$ est 
surjectif;
de plus, on a $f_{*}(\O_{S_{\epsilon}})=\O_{\M_c}$.
\end{theor}

{\bf Preuve:}
Si $F$ est un faisceau stable de classe $c$, la condition
$n<(l+1)(l+2)$ signifie que $\chi(F(l))>0$ et par suite on peut
considérer les systèmes cohérents $\Lambda=(\Gamma,F(l))$ avec
$\Gamma\subset \H^0(F(l))$. Ces systèmes cohérents sont 
obligatoirement
$\epsilon$--stables. Il en résulte que la fibre de $f$ au-dessus du
point défini par $F$ est isomorphe à l'espace projectif (des 
droites)
$\proj(\H^0(F(l)))$. Le résultat est alors évident si $n$ est impair, car
il n'y a alors que des points stables: si ${\mathcal F}$ est une famille
universelle, l'espace de modules s'identifie d'ailleurs dans ce cas au
schéma de Grothendieck 
$$\proj(\underline{Ext}^2_{pr_1}({\mathcal F},pr_2^*(\omega_{\pp})))$$
où $\omega_{\pp}$ désigne le fibré canonique de $\pp$
($\omega_{\pp}\simeq\O(-3)$).

Dans le cas où $n$ est pair, on peut employer l'argument suivant :
l'espace de modules $S_{\epsilon}$ étant  intègre, on sait que l'algèbre $\A=f_*(\O_{S_{\epsilon}})$ est intègre.
Mais alors $\M'=\spec(\A)$ est une variété projective intègre au-dessus de $\M_c$, et le \mor $\M'\to\M_c$ est projectif et birationnel. De plus, la variété $\M_c$ est normale: ceci implique que l'image directe du \fs structural de $\M'$ est $\O_{\M_c}$ (\cite{Har}, chap. III, cor. 11.4). Le \th en résulte. $\Box$

\subsection{Le fibré déterminant}

Si $a$ n'est pas valeur critique, on sait qu'il existe un système
cohérent universel paramétré par l'espace de modules $S_a$. On l'écrit
sous la forme $({\mathcal V},{\mathcal F}(l))$, où ${\mathcal V}$ est un fibré inversible sur $S_a$, et ${\mathcal F}$ une famille plate de faisceaux cohérents de classe $c$ paramétrée par $S_a$. On peut donc, pour toute classe
$u\in c^{\bot}$ dans $K(\pp)$, de dimension 1, (l'orthogonal est
pris relativement à la forme quadratique  sur $K(\pp)$) définir un fibré déterminant $\D_{a,u}$ sur $S_a$ par la formule
$$\D_{a,u}=\lambda_{{\mathcal F}}(-u)=\det\, p_{1!}({\mathcal F}\cdot p_2^*(-u)).$$
Dans cette formule $p_1$ et $p_2$ sont les projections canoniques:
$$\begin{array}{ccc}
 {S_a\times\pp} & \rightto{p_2} & {\pp} \\
 \downto{p_1} & & \\
 {S_a}     &  & \\
\end{array}$$
et ${\mathcal F}(u)={\mathcal F}\cdot p_2^*(u)$  désigne
la classe (dans le groupe de Grothendieck
$K(S_a\times \pp)$ des classes de faisceaux algébriques
cohérents plats sur $S_a$) produit de la classe de ${\mathcal F}$  par l'image réciproque de $u$  par la projection $p_{2}.$  Le morphisme 
$$p_{1!}: K(S_a\times \pp)\ra K(S_a)$$ 
 est le morphisme  qui
associe à la classe d'un faisceau  algébrique cohérent ${\cal
F}$ plat sur $S_a$
 la classe de  \hbox{$\sum_{q}(-1)^{q}R^qp_{1*}({\mathcal F}).$} Ces \fx de
 cohomologie sont les \fx de cohomologie d'un complexe fini de fibrés
 vectoriels $Rp_{1*}({\mathcal F})$.
Par la propriété universelle du fibré déterminant, on a
$f^*(\D_u)=\D_{{\epsilon},u}$. Désignons par $\mathfrak{d}$ le fibré déterminant sur $\Hilb^{n+l^2}(\pp)$, et considérons l'ouvert $\U$ de
ce schéma de Hilbert où le faisceau ${\mathcal R}$ est localement
libre. Cet ouvert est l'espace tout entier pour $l=1$ et
  il a son complémentaire de codimension
$2l$, pour $l$ supérieur à $1$. Il est invariant sous l'action du groupe $\sl(3)$.

\begin{theor}
\label{th1}
Soit $n$ un entier $\ge 3$. Soit $l$ un entier $>0$ tel que $l(l-1)\le
n<(l+1)(l+2)$. Alors on a un isomorphisme de $\sl(3)$--représentations
$$\H^0(\M_c,\D)=\H^0(\U,\es^l{\mathcal R}\otimes{\mathfrak{d}}).$$
\end{theor}

{\bf Preuve:}
D'après le théorème \ref{epsi}, on a
$$\H^0(\M_c,\D)=\H^0(S_{\epsilon},\D_{\epsilon,\mathfrak{u}}).$$
Ici, $\mathfrak{u}$ désigne la classe du faisceau $\O_H(-1)$, 
orthogonale à $c$, où $H$ est une droite de $\pp$. Maintenant, d'après
le résultat de Min He, les espaces de modules $S_a$ sont des 
variétés
normales: les espaces vectoriels de sections restent inchangés par
restriction à un ouvert dont le complémentaire est de codimension 
$\ge
2$. Puisque les fermés $\Sigma,\Sigma_{-}$, et $\Sigma_{+}$ du
théorème \ref{minhe} sont de codimension $\ge 2$ on voit que la
représentation $\H^0(S_a,\D_{a,\mathfrak{u}})$ est indépendante de 
$a$. Il reste
à voir ce qu'est cette représentation pour $a=a_{max_{-}}$. Ceci
résulte du calcul du fibré déterminant $\D_{a_{max_{-}},{\mathfrak{u}}}$ sur
l'espace de modules $S_{a_{max_{-}}}$.

\begin{lemme}
Soit ${\mathfrak{a}}$ le fibré tautologique $\O_{\proj({\mathcal 
R})}(1)$
sur $\proj({\mathcal R})$. Alors
$$\D_{a_{max_{-}},\mathfrak{u}}={\mathfrak{a}}^{\tens 
l}\tens\pi^*(\mathfrak{d}).$$
\end{lemme}

\par{\bf Preuve du lemme~:}
Rappellons que $\Xi\subset\Hilb^{n+l^2}(\pp)\times\pp$ est le
sous-schéma universel, ${\mathcal I}$ le faisceaux d'idéaux associé et
${\mathcal R}$ est le 
faisceau
algébrique cohérent $R^1pr_{1*}({\mathcal I}(2l-3))$. Considérons
l'extension canonique sur $\proj({\mathcal R})\times \pp$~:
$$0\to{\mathfrak{a}}\boxtimes\O(-l)\to\maF\to(\pi\times id_{\pp})^*(\I(0,l))\to
0$$
où $\pi:\proj(\R)\to \Hilb^{n+l^2}(\pp)$ est le \mor canonique.
La classe ${\mathfrak{u}}$ est celle du \fs $\O_H(-1)$, où $H$ est une
droite de $\pp$. On a alors par changement de base
\begin{eqnarray*}
 \D_{a_{max_{-}},{\mathfrak{u}}}&=&\O(l{\mathfrak{a}})\tens\pi^*(\lambda_{\O_{\Xi}}({\mathfrak{u}}(l))\\
&=&\O(l{\mathfrak{a}})\tens\pi^*(\lambda_{\O_{\Xi}}({\mathfrak{u}}))
\end{eqnarray*}
d'après \cite{LeP-Durham}, prop. 2.9. Le \fib inversible
$\lambda_{\O_{\Xi}}({\mathfrak{u}})$ est le \fib déterminant de Donaldson
$\de$. Ceci démontre le lemme. $\Box$

Il reste maintenant à enlever le fermé image réciproque du lieu de
Brill-Noether $B$, qui est de codimension $\ge 2$, pour obtenir le
résultat du théorème.$\Box$


\section{Sections de $\es^l{\mathcal R}\tens\de$  sur le schéma de 
Hilbert}
\label{resolution}

Le but de cette section est de terminer la démonstration du 
théorème
\ref{teore}.
On supposera partout dans la suite que $n\ge 3$.
On sait par le théorème \ref{th1} que pour $l(l-1)\le n<(l+1)(l+2)$, 
on a un isomorphisme $\H^0(\M,\D)=\H^0(U,\es^l\R\tens\de)$ où $U$ 
est l'ouvert du schéma de Hilbert $\Hilb^{n+l^2}(\pp)$ où le \fs 
$\R$ est localement libre. Cet ouvert a son complémentaire de 
codimension $\ge 2l$. Il est invariant sous l'action de ${\sl(3)}$.

Lorsqu'il n'est pas spécifié, les produits tensoriels de \fx \algs
considérés sont des produits tensoriels sur le \fs structural du
schéma de base. Les \fibs vectoriels sont identifiés à des \fx
localement libres de rang fini.

On va se concentrer maintenant sur ce nouvel espace de sections  
$\H^0(U,\es^l\R\tens\de)$. On note $m=n+l^2$.

Soit $\Xi\subset\hil\times\pp$ le  sous-schéma universel, $pr_1:\Xi\to\hil$ et $pr_2:\Xi\to\pp$ les deux projections, et 
$\I_{\Xi}$ le \fs d'idéaux associé. Le \fs $\R$ est défini comme 
$\R=R^1pr_{1*}(\I_{\Xi}(2l-3))$. On note $k=2l-3$. En partant de la 
suite exacte fondamentale associée à $\Xi$ sur  
$\Hilb^m(\pp)\times\pp$:
$$0\ra\I_{\Xi}\ra\O\ra \O_{\Xi}\ra 0$$
tensorisée par $pr_2^*(\O_{\pp}(k))$ et restreinte à l'ouvert 
$U\times\pp$ de $\Hilb^m(\pp)\times\pp$  on trouve la suite exacte
$$0\ra\I_{\Xi}(k)|_{U\times\pp}\ra\O(k)|_{U\times\pp}\ra\O_{\Xi}(k)|_{U\times
\pp}\ra 0$$
Par image directe sur $U$ par la projection $pr_{1}$, on obtient une présentation $\sl(3)$-équivariante 
de $\R$ sur $U$:
\begin{equation*}
0\ra pr_{1*}(\I_{\Xi}(k))\ra pr_{1*}(\O(k))\ra pr_{1*}(\O_{\Xi}(k))\ra \R\ra R^1
pr_{1*}(\O(k))\ra 0
\end{equation*}
 
Par définition de $U$ on a  $\H^0(\pp,\I_Z(k))=0$ pour les schémas
 $Z\in U$,
 donc par le \th de 
semi-continuité on a une suite exacte de \fx localement libres sur $U$:
\begin{equation}
0\ra \H^0(\pp,\O(k))\tens \O_{U}\ra
pr_{1*}(\O_{\Xi}(k))|_{U}\ra\R\ra 0
\label{*}
\end{equation}

On obtient par suite une résolution $\sl(3)$-équivariante de
$\es^l\R\tens\de$ par le complexe  de Koszul $K^{\cdot}$ défini en degré $i$ par $K^{
i}=\Lambda^{-i}\es^kE\tens\es^{l+i}\vkde$ pour $i=0,\ldots,l$ où 
$\vk=pr_{1*}(\O_{\Xi}(k))$, $E=\H^0(\pp,\O(1))$~:
\begin{eqnarray}
\label{**}
 K^{\cdot}\to \es^l\R\tens\de.
\end{eqnarray}
Par suite la cohomologie de 
$\es^l\R\tens\de$  se calcule à l'aide de la suite spectrale
$E^{p,q}_2=\H^q(K^{-p})$  dont 
l'aboutissement en degré $0$ est $\H^0(U,\es^l\R\tens\de)$.

Pour les \fx localement libres $K^{i}$ on va se placer indifféremment sur les restrictions à $U$ ou sur $\hil$ tout entier, 
puisqu'on s'intéresse seulement à la cohomologie de ces  \fx  jusqu'en degré $l$; 
comme le complémentaire de $U$ est de codimension $\ge 2l$, celle-ci coïncide sur $U$ 
avec la cohomologie sur tout $\hil$ jusqu'en degré $2l-2$, par les propriétés de la cohomologie locale (\cite{Hart}). Pour $l=1$, $U=\hil$.

\begin{lemme}
\label{Kawamata}
Pour $q>0$, $\H^q(\Hilb^m(\pp),\de)=0$ et $\H^0(\hil,\de)=\es^mE$.
\end{lemme}

\par{\bf Preuve du lemme~:}
Le résultat découle par le \th de Kawamata-Viehweg (voir \cite{C-K-M}, page 52, thm.8.3.), $\de$ 
étant le fibré déterminant, big et nef, sur le schéma de Hilbert vu comme espace de modules de 
\fx sans torsion de rang 1 sur $\pp$, de classes de Chern $c_1=0, c_2=m$.

Soit $\pi$ le morphisme de Hilbert-Chow, $\pi:\Hilb^m(\pp)\to
\es^m(\pp)$ qui associe à un schéma fini $Z$ le cycle
$\sum_{x\in\pp}lg\,Z_xx$  dans $\es^m(\pp)$, la puissance symétrique $m$-ième de
$\pp$. Naturellement, l'espace $\es^m(\pp)$ est le quotient de la
puissance $m$-ième $\pp^m$ de $\pp$ par le groupe symétrique
$\sigm$. Il est constitué de cycles, combinaisons linéaires
de points distincts $x_i$ de $\pp$, 
$$\sum_{\sum_i  \lambda_i=m}\lambda_i[x_i],$$
à coefficients $\lambda_i>0$.
 Le support d'un schéma fini de longueur
$m$ est un tel cycle, si on tient compte des multiplicités des
points. Le morphisme de Hilbert-Chow vérifie
 $\pi_*\O_{\hil}=\O_{\es^m(\pp)}$ et (cf. \cite{LeP-Durham})
\begin{equation}
\label{delta}
\pi^*(\O(1,1,\ldots,1)^{\sigm})=\de.
\end{equation}
 On obtient alors $\H^0(\Hilb^m(\pp),\de)=\H^0(\es^m(\pp),\pi_*(\de))=\H^0(\pp^m,\O(1,\cdots,1))^{\sigm}=\es^m\H^0(\pp,\O(1))=\es^mE$.$\Box$

\vs Dans \cite{D}, il est démontré que
$\H^0(\hil,\vkde)=\es^{2l-2}E\tens\es^{m-1}E$ et que $\H^q(\hil,\vkde)=0$
pour $q>0$. Comme annoncé dans l'introduction (th. \ref{enunt}),  on
montrera aussi que 
$\H^1(\hil,\es^2\V_k\tens\de)=0$ et on calculera 
$\H^0(\hil,\es^2\V_k\tens\de)$ et $\H^0(\hil,\es^3\V_3\tens\de)$. On  peut 
ainsi calculer $\H^0(\hil,\es^l\R\tens\de)$ pour $l=1,2,3$. Pour aller 
plus loin on se heurte à des difficultés liées au calcul des 
\hbox{$\H^q(\hil,\es^l\vkde)$} pour $q=0,  l>3$ ou $q>0, l>1$.
\vs Ceci  limite le  calcul du $\H^0(U,\es^l\R\tens\de)$ à $l=3$ 
ce qui restreint les valeurs de $n$ à $n\le 19$.
\vs On commence par remarquer que pour le calcul d'un espace de sections 
d'un fibré sur $\hil$ on peut se placer sur un grand ouvert de $\hil$, 
pourvu que cet ouvert ait un complémentaire de codimension au moins 2. 
Ceci est le cas pour l'ouvert $\hils$ formé par les schémas avec au 
plus un point multiple, qui soit double, soit les schémas dont le 
cycle correspondant 
est $x_1+x_2+\cdots+x_m$ ou $2x_1+x_3+x_4+\cdots+x_m$ avec $x_i$ distincts.
On note $\es^m_*(\pp)$ l'ouvert des cycles de cette forme.
 L'avantage d'utiliser $\hils$ 
est qu'on peut le décrire comme quotient $q$ 
de l'éclaté $B$ de $\ppms=p^{-1}(\sms)$ (où $p:\ppm\to\sm$ est le 
quotient  de $\ppm$ sous l'action du groupe symétrique ${\mathfrak{S}}_m$ ) selon la réunion $D$ 
des diagonales $\Delta_{ij}=\{(x_1,\cdots,x_m)\in\ppms|x_i=x_j\}$ pour 
$i<j$, disjointes dans $\ppms$. On note $\rho$ cet éclatement. On a un 
diagramme commutatif:
$$\begin{array}{ccc}
{B} & \rightto{\rho} & {\ppms} \\
\downto{q} & & \downto{p} \\
{\hils} & \rightto{\pi} &{\sms}
\end{array}$$

On  montrera comment, à l'aide de cette description, on peut ramener 
les calculs de la \coh des fibrés sur $\hils$ à des calculs des 
invariants de la \coh de certains \fx sur $\ppms$. On utilise les mêmes notations $\hilxs$ et $\xms$ pour une surface quelconque $X$.


\subsection{Une filtration}

\label{filtration}

On introduit ici des notations et des résultats très utiles pour la 
suite. On va se placer dans le cadre général d'une variété \alg 
lisse $M$, munie d'un fibré de rang r, $W$, et d'un \fib  $L$ 
sur une sous-variété lisse $D$ de $M$. On note $W_D$ la restriction de $W$ à $D$.

On considère un \mor \sur $\epsilon:W\surto L$. Le noyau de ce \mor 
définit un \fs sans torsion $V$. Ce \mor induit un \mor \sur de \fibs 
en algèbres graduées $\sym\epsilon:\sym W=\oplus_{i\ge 
0}\es^iW\surto\sym L=\oplus_{i\ge 0}\es^iL$, noté encore 
$\epsilon$. On note $I$ le faisceau noyau.

Considérons la filtration $\F^k\sym W=I^k\sym W$, pour $k\ge 0$; elle 
est compatible avec la graduation.

\begin{prop}
\label{3.1}
Soit  $\N$ le \fib conormal de $D$ dans $M$ et $K$ le noyau du \mor 
canonique $\epsilon|_D$ (noté encore $\epsilon$)$:W_D\to L$.
Le gradué associé à cette filtration est donné par
\begin{equation}
\gr_q(I^p/I^{p+1})=\es^q\N\tens\es^{p-q}K[-p+q]\tens\sym L
\label{filtratio}
\end{equation}
si $p\ge q\ge 0$ et $0$ sinon.
\end{prop}

\par{\bf Preuve:}

Afin de décrire le gradué associé à cette \fil on va considérer $\sym W$ comme image directe de l'algèbre des 
fonctions régulières $\O_{W^*}$ sur l'espace total du \fib dual 
$W^*$, par la projection canonique $p:W^*\to M$, et $\sym L$  comme 
image directe de l'algèbre des fonctions régulières $\O_{L^*}$ sur 
l'espace total du \fib dual $L^*$, par la restriction de $p$ à 
$L^*\subset W^*$, notée encore $p:L^*\to D$.

On désigne par $W_D$  la restriction de $W$ à $D$. Considérons les 
inclusions de variétés lisses
\begin{equation}
L^*\subset W^*_D\subset W^* 
\label{s1}
\end{equation}
et désignons par $\I$ l'idéal de $L^*$ dans $W^*$, par $\J$ 
l'idéal  de $W^*_D$ dans $W^*$, et par $\ID$ l'idéal de $L^*$ dans 
$W^*_D$. On a une suite exacte $0\to\J\to\I\stackrel{r}{\to} \ID\to 0$ 
où $r$ est le \mor de restriction. On considère les \fils de 
$\O_{W^*}$ et $\O_{W^*_D}$ définies par les puissances des idéaux  $\I$ et $\ID$.

Comme $p$ est un \mor affine, $W^*$ est un schéma affine sur $M$, il 
y a une correspondance entre les \fx d'idéaux de $\O_{W^*}$ et les 
idéaux de $\sym W$ donnée par $\I\mapsto p_*(\I)$. L'idéal $\I$ se 
correspond ainsi à $I$ et $\I^k$ à $I^k$. De cette fa\c con, la \fil 
de $\sym W=p_*(\O_{W^*})$ définie par image directe coïncide 
avec la \fil définie par l'idéal $I=p_*(\I)$ noyau du \mor 
$\epsilon$.
Les \fibs conormaux correspondants aux inclusions (\ref{s1})
s'écrivent 
\hbox{$\N_{L^*/W^*}=N^*_{L^*/W^*}=\I/{\I^2}$}, 
\hbox{$\N_{W^*_D/W^*}=N^*_{W^*_D/W^*}=\J/{\J^2}$}, 
\hbox{$\N_{L^*/W^*_D}=N^*_{L^*/W^*_D}=\ID/{\ID^2}$} et on  a une suite 
exacte
\begin{equation}
0\to \N_{W^*_D/W^*}|_{L^*}\to \N_{L^*/W^*}\to \N_{L^*/W^*_D}\to 0 
\label{s2}
\end{equation}

Le fibré conormal $\N_{W^*_D/W^*}|_{L^*}$ s'identifie à $p^*(\N)$ et $\N_{L^*/W^*_D}$ 
à $p^*(K)$. La suite (\ref{s2}) devient
\begin{equation}
0\to p^*(\N)\to \I/{\I^2}\to p^*(K)\to 0.
\label{s3}
\end{equation}
À partir de cette suite exacte on obtient une \fil décroissante du 
$\O_{L^*}$-module 
$\es^p(\I/{\I^2})=\I^p/{\I^{p+1}}=\gr_p(\O_{W^*})$
 par des $\O_{L^*}$-modules
$$\F^q(\es^p\N_{L^*/W^*}) \\=\\ \Im((\J/{J^2})^{\tens q}\tens(\I/{\I^2})^{\tens(p-
q)}\to\es^p(\I/{\I^2})) \\=\\ \Im(\J^q\I^{p-q}\to\I^p/{\I^{p+1}})$$
 si $p\ge q\ge 0$ et $0$ sinon, de gradué associé
$$\gr_q(\es^p(\I/{\I^2}))=\es^q\N_{W^*_D/W^*}\tens\es^{p-q}\N_{L^*/W^*_D}=p^*(\es^q\N\tens\es^{p-q}K)$$
si $p\ge q\ge 0$ et $0$ sinon.
Par application du foncteur image directe $p_*$ qui est exact puisqu'il 
s'agit d'un morphisme affine, on obtient une \fil de 
$p_*(\I/{\I^{p+1}})=I^p/{I^{p+1}}$ dont le gradué en degré $q$ 
est le $\sym L$-module gradué fourni par la formule de projection
$$\gr_q(I^p/{I^{p+1}})=\es^q\N\tens\es^{p-q}K[-p+q]\tens\sym L$$
si $p\ge q\ge 0$ et $0$ sinon.
Pour comprendre le décalage qui apparaît dans la graduation de 
$p_*(p^*(\es^q\N\tens\es^{p-q}K))=\es^q\N\tens\es^{p-q}K\tens\sym L$ il 
faut comprendre l'action de $\comp^*$ sur le \fib conormal 
$\N_{L^*/W^*_D}=p^*(K)$ et sur le \fib conormal $\N$. L'action  de 
$\comp^*$ sur $L^*$, $W^*_D$ et $W^*$ est par homothétie, d'où une 
action sur les trois \fibs normaux et respectivement conormaux .
Sur $M$ et $D$, et par conséquent  sur $\N$, $\comp^*$ agit 
trivialement. Donc la composante homogène de degré $i$ de 
$p_*(p^*(\N))$ est $\N\tens\sym^iL$.
Sur $\N_{L^*/W^*_D}=p^*(K)$, l'action est donnée par 
$\lambda\cdot(x,v)=(\lambda x,\lambda^{-1}v)$ pour $x\in L^*$ et $v\in 
K_{p(x)}$. Donc la composante homogène de degré $i$ de 
$p_*(p^*(\es^pK))$ est donnée par $\es^pK\tens\sym^{i-p}L$ si $i\ge p$ $\Box$.

\begin{lemme}
\label{l2}
L'image de $\es^kV$ dans $\sym W$ par le \mor $\es^k\iota$ engendre comme idéal $I^k$ en degré $\ge k$  où $\iota$ est l'inclusion de $V$ dans 
$W$.
\end{lemme}
 
\par{\bf Preuve: } 
On va montrer plus précisément que pour $i\ge k$, l'image du \mor naturel
  $v_{i,k}:\es^kV\tens\es^{i-k}W\to\es^iW$ est exactement  $(I^k)_i=I^k\bigcap\es^iW$.
On commence par un
\begin{slemme}
\label{souslemme}
Si $n\ge 1$, l'image du \mor $v:V\tens\es^{n-1}W\to\es^nW$  (qui à $v\tens w_1\cdots w_{n-1}$ 
associe le symétrisé de $\iota(v),w_1,\ldots,w_{n-1}$ dans $\es^nW$, noté $\iota(v)w_1\cdots w_{n-1}$ 
) est exactement $I_n=\Ker(\es^n\epsilon:\es^nW\to\es^nL)$.
\end{slemme}
 
\par{\bf Preuve du sous-lemme:  }
En tensorisant $n$ fois avec elle-même la suite exacte sur $M$:
$$ V\stackrel{\iota}{\to}W\stackrel{\epsilon}{\to}L\to 0$$ 
on trouve une suite exacte
$$T=V\tens W^{\tens(n-1)}\oplus W\tens V\tens
W^{\tens(n-2)}\oplus\cdots\oplus W^{\tens(n-1)}\tens
V\stackrel{}{\to}W^{\tens n}\stackrel{\epsilon^{\tens n}}\to L^{\tens
  n}\to 0.$$
Le cas $n=2$ est la proposition 6, page 8, chap. III-1 de \cite{Bourbaki}. Ce cas se généralise sans peine par
récurrence au produit tensoriel d'un nombre fini quelconque de suites.
Il suffit alors de remarquer que la suite exacte des invariants par le
groupe symétrique ${\mathfrak{S}}_n$ reste exacte et on a exactement que
$$T^{{\mathfrak{S}}_n}=V\tens\es^{n-1}W\stackrel{v}{\to}\es^nW=(W^{\tens
n})^{{\mathfrak{S}}_n}\stackrel{\es^n\epsilon}{\surto}\es^nL=(L^{\tens
n})^{{\mathfrak{S}}_n}\to 0$$
est exacte d'où $I_n=\Ker\es^n\epsilon=\Im v$. $\Box$

\par{\bf Preuve du lemme \ref{l2}: }

Une section locale $\alpha$ de $(I^k)_i$ s'écrit comme 
$\alpha=\sum\alpha_1\cdots\alpha_k$  avec $\alpha_j\in I_{n_j}$ et 
$\sum^k_{j=1}n_j=i$. Par le sous-lemme \ref{souslemme} , chaque 
$\alpha_j$ avec $n_j\ne 0$ provient d'une section $\widetilde{\alpha_j}$ de 
$V\tens \es^{n_j-1}W$, et si $n_j=0$, $\alpha_j$ est une section de 
$I_D=I_0$. Mais $\epsilon(I_D\tens W)=0$ donc $I_DW\subset 
\iota(V)$. Au total, $\alpha$ provient d'une section de 
$\es^kV\tens\es^{i-k}W$.$\Box$

\begin{rem}
{\rm Considérons maintenant un \fib inversible $\A$ sur $M$. Alors on a un \mor de $\sym W$-modules gradués
$${\mathcal{M}}=\sym W\tens\A\to\sym L\tens\A$$
de noyau $I{\mathcal{M}}$. Considérons la \fil $I^k{\mathcal{M}}$. Cette \fil est compatible avec la graduation et le \mor canonique
$$\Phi:\es^kV|_D\tens_{\O_D}\sym L[-k]\tens\A|_D\to I^k{\mathcal{M}}/I^{k+1}{\mathcal{M}}$$
est un \iso en degré $\ge k$.}
\end{rem}


\subsection{Éclatement de $M$ le long de $D$}

On considère l'éclatement $\rho:\widetilde{M}\to M$ de $M$ le long de 
$D$, et les images réciproques $\widetilde{W}$ et $\tL$ de $W$ et $L$ par 
$\rho$: $\tW=\rho^*(W), \tL=\rho^*(L).$ On note $\tV$ le noyau du 
\mor \sur, noté encore $\epsilon$, de $\tW$ dans $\tL$. Puisque le 
support de $\tL$ est un diviseur (le diviseur exceptionnel ${\bf E}$), $\tV$ 
est localement libre. De manière analogue, on considère le 
noyau $\tI$ de $\epsilon:\sym\tW\to\sym\tL$, et la \fil $\tI^k$ de $\sym\tW$.

\begin{lemme}
\label{l3}

\begin{itemize}

\item{\rm (i)} Le \mor canonique $\rho^*:\sym W\to\rho_*(\sym \tW)$ induit un \iso $I^k\stackrel{\sim}{\to}\rho_*(\tI^k)$
\item{\rm (ii)} Les images directes $R^q\rho_*(\tI^k)$ sont nulles pour $q>0$.
\end{itemize}
\end{lemme}

\par{\bf Preuve:}

L'éclatement $\rho$ vérifie $\rho_*(\O_{\tM})=\O_M$
 et $R^q\rho_*(\O_{\tM})=0$ pour $q>0$, d'après le lemme 3.5 de
 \cite{SGA}, exposé VII.
On a alors, par la formule de projection, un \mor
$$\sym W\stackrel{\sim}{\to}\rho_*(\rho^*\sym W)=\rho_*(\sym\tW)$$
qui 
est un isomorphisme et de même pour $\sym L\stackrel{\sim}{\to}\rho_*(\sym\tL)$. 
On a donc un diagramme commutatif
$$\begin{array}{ccccccccc}   
0 & \rightto{} & {I} &  \rightto{} & {\sym W} & \rightto{\epsilon} & {\sym L} & \rightto{}
& {0} \\
& & & &\downto{\sim}&  &\downto{\sim} & & \\
{0} & \rightto{} & {\rho_*(\tI)} & \rightto{} 
&{\rho_*(\sym\tW)} &\rightto{\rho_*\epsilon}
& {\rho_*(\sym\tL)} & \rightto{} & {\ldots}
\end{array}$$
qui nous assure que $\rho_*\epsilon$ est \sur et qu'il y a aussi un 
isomorphisme $I\stackrel{\sim}{\to}\rho_*(\tI)$. D'où un \mor 
$I^k\to\rho_*(\tI^k)$.

On suppose par récurrence que pour tout $i\le k$ on a le résultat 
(pour $k=0$ ceci est clair: $\sym W\stackrel{\sim}{\to}\rho_*(\sym\tW)$
et $R^q\rho_*(\sym\tW)=0$ pour $q>0$) et on va le prouver pour $k+1$. 
On a un \mor de suites exactes
$$\begin{array}{ccccccccccc}
{0}&\rightto{}&{I^{k+1}}&\rightto{}&{I^k}&\rightto{}&{I^k/{I^{k+1}}}&\rightto{}&{0}& &\\
&&\downto{a}&&\downto{\sim,b}&&\downto{c}&&&&\\
{0}&\rightto{}&{\rho_*(\tI^{k+1})}&\rightto{}&{\rho_*(\tI^k)}&\rightto{d}&{\rho_*(\tI^k/{\tI^{k+1})}}&\rightto{}&{R^1\rho_*(\tI^{k+1})}&\rightto{}&0\
\end{array}$$
où $b$ est un isomorphisme. On commence par prouver que $c$ est un 
isomorphisme et que \linebreak $R^q\rho_*(\tI^k/{\tI^{k+1}})=0$ pour $q>0$. On 
en déduira que $d$ est \sur d'où $a$ sera un isomorphisme et 
$R^q\rho_*(\tI^{k+1})=0$ pour $q>0$.

\vs Mais on a construit dans  la  précédente section une filtration de 
chacun des \fx $I^k/{I^{k+1}}$ et $\tI^k/{\tI^{k+1}}$ de gradués connus. 
Ces \fx sont supportés par $D$ et ${\bf E}$ respectivement et le \mor $\pi$ 
en restriction à $D$ s'écrit comme $\rho:{\bf E}={\Bbb P}(\N^*_D)\to D$. 
Le noyau $\tK$ de 
$\tW|_{{\bf E}}\to\tL$ s'identifie à l'image réciproque de $K$.  Le \fib conormal à ${\bf E}$, $\N_{{\bf E}}$, est dans ce cas le \fib 
$\O(1)$ relatif sur cet espace projectif (on a pris le projectif de 
Grothendieck). La \fil $\tF^j$ de $\tI^k/{\tI^{k+1}}$ est de gradué
$$\gr_j(\tI^k/{\tI^{k+1}})=\es^j\N_{{\bf E}}\tens \es^{k-j}K[-k+j]\tens\sym\tL$$ 
si $k\ge j\ge 0$ et $0$ sinon.

Comme
$R^q\rho_*(\gr_j(\tI^k/{\tI^{k+1}}))=R^q\rho_*(\es^j\N_{{\bf E}})\tens\es^{k-j}K
[-k+j]\tens\sym L$
et que $\rho_*(\es^j\N_{{\bf E}})=\es^j\N_D$ et $ R^q\rho_*(\es^j\N_{{\bf E}})=0$ pour 
$q>0$ on obtient une \fil  $\F^j=\rho_*(\tF^j)$ de 
$\rho_*(\tI^k/{\tI^{k+1}})$ de gradué
$$\rho_*(\gr_j(\tI^k/{\tI^{k+1}}))=\es^j\N_D\tens\es^{k-j}K[-k+j]\tens\sym L=\gr_j(I^k/{I^{k+1}})$$ 
si $k\ge j\ge 0$ et $0$ sinon, et telle que 
$R^q\rho_*(\F^j)=0$ si $q>0$, pour tout $j$. En particulier pour 
$j=0$ on obtient que $R^q\rho_*(\tI^k/{\tI^{k+1}})=0$ pour $q>0$. 
Le \mor $c$ est compatible avec les \fils et induit l'identité sur les 
gradués, d'où aussi l'isomorphisme
$c:I^k/{I^{k+1}}\stackrel{\sim}{\to}\rho_*(\tI^k/{\tI^{k+1}})$.

\begin{cor}
\begin{itemize}
\item{\rm (i)} L'image de l'inclusion canonique 
$\phi:\rho_*(\es^k\tV)\hookrightarrow\es^k W$ est exactement $(I^k)_k$.
\item{\rm (ii)} $R^q\rho_*(\es^k\tV)=0$ pour $q>0$.
\end{itemize}
\end{cor}

\par{\bf Preuve:}

Cela revient à écrire les résultats du lemme \ref{l3} en degré 
$k$ en tenant compte du lemme \ref{l2} et du fait que 
$\es^k\widetilde{\iota}:\es^k\tV\hookrightarrow\es^k\tW$ reste une inclusion, ou $\widetilde{\iota}$ est l'inclusion de $\tV$ dans $\tW$.$\Box$


\subsection{Calculs de cohomologie sur $\hils$}
\label{calco}

On appliquera ici les résultats des deux sections précédentes à 
notre situation particulière. On n'aura pas besoin ici de se placer 
sur le plan projectif. Les résultats restent valables sur une surface 
\alg lisse quasi-projective quelconque $X$. La description de $\hilxs$ se 
fait alors exactement comme pour $\pp$, en utilisant l'éclatement $B$ 
de $\xms$:
$$\begin{array}{ccc}
{B}&\rightto{\rho}&{\xms}\\
\downto{q}& &\downto{p}\\
{\hilxs} & \rightto{\pi} & {\smxs}
\end{array}$$

On considère plus généralement le \fib $\V_L$ sur $\hilx$ 
associé à un \fib  $L$ sur $X$,\linebreak 
$\V_L=pr_{1*}(\O_\Xi\tens pr^*_{2}(L))$ où $pr_1$ et $pr_2$ sont les 
deux projections du schéma universel $\Xi\subset\hilx\times X$ sur 
$\hilx$ et respectivement $X$.

\vs On garde les notations introduites juste avant la section \ref{filtration} pour les diagonales $D$ et $\Delta_{ij}$ 
de $\xms$. Sur $B$, le diviseur exceptionnel ${\bf E}$ se décompose en 
composantes disjointes ${\bf E}=\bigcup_{i<j}{\bf E}_{i,j}$. Alors le schéma 
universel $\Xi_B\subset B\times X$, paramétré par $B$, a $m$ 
composantes irréductibles $\Xi_i$ et la projection 
$pr_1:\Xi_i\bigcap\Xi_j\to {\bf E}_{i,j}$ est un isomorphisme. On en déduit 
une suite exacte sur $B\times X$:
\begin{equation}
0\to\O_{\Xi_B}\to\oplus_i\O_{\Xi_i}\to\oplus_{i<j}\O_{{\bf E}_{i,j}}\to 0 
\label{s6}
\end{equation}
et comme, par changement de base, $q^*(\V_L)=pr_{1*}(\O_{\Xi_B}\tens 
pr_2^*(L))$, on a, après tensorisation par 
$pr^*_2(L)$ de la suite (\ref{s6}) et image directe par $pr_1$, une suite exacte sur $B$:
$$0\to q^*(\V_L)\to\oplus_ip_i^*(L)\to\oplus_{i<j}p^*_{i,j}(L_\Delta)\to 
0$$
où $p_i$ désigne aussi bien la $i$-ème projection $\xms\to X$ que 
sa composée avec $\rho:B\to X$; de même pour $p_{i,j}:\xms\to 
X\times X$ et $p_{i,j}:B\to X\times X$. 

Le sous-schéma $\Xi_i$ est 
l'image réciproque de la diagonale $\Delta$ de $X\times X$ par 
l'application $(p_i,id_X)$. Le \fib $L_\Delta$ est l'image réciproque 
de $L$ par l'une des projections de la diagonale de $X\times X$ sur $X$, 
qui sont des isomorphismes.

Le \fib $\tmL=\oplus_{i<j}p^*_{i,j}(L_\Delta)$ a pour support le 
diviseur exceptionnel ${\bf E}$. Il est l'image réciproque  par $\rho$ du 
\fib $\L=\oplus_{i<j}p^*_{i,j}(L_\Delta)$ sur $\xms$, dont le support est 
 $D$.

On note aussi $\tW$ et $W$ le \fib $\oplus_ip^*_i(L)$ sur $B$ et 
 sur $\xms$ respectivement. On reconnaît maintenant la variété 
$M=\xms$ et $\tM=B$, et le \fib $\tV=q^*(\V_L)$. Le groupe symétrique
$G=\sigm$ opère sur la situation. Toutes les \fils qui
interviennent sont invariantes sous l'action de $G$, et les morphismes sont 
$G$-équivariants. Comme le \mor $\epsilon$ est donné par 
$(s_i)_i\mapsto(s_i|_D-s_j|_D)_{i,j}$, l'action induite sur $\L$ et 
$\tL$ est telle que la transposition $\tau_{i,j}$ qui échange $i$ et 
$j$ change le terme d'indice $(i,j)$ en son opposé. Le groupe $G$ 
étant fini, la \coh du \fs des invariants $F^G$ sur $\hilxs$ (ou sur 
$\smxs$), où $F$ est un $G$-\fs \alg cohérent sur $B$ (ou sur $\xms$) s'identifie à la 
\coh équivariante de $F$, c'est-à-dire aux invariants de la \coh de $F$.

\vs Le \th suivant nous montre comment on peut ramener le calcul de la \coh 
de $\es^l\V_L$ sur $\hilxs$ à un calcul de \coh équivariante sur 
$\xms$:

\begin{theor}

\begin{itemize}

\item{\rm (i)} Il existe une inclusion canonique 
$\pi_*(\es^l\V_L)\hookrightarrow(\sym W)^G$ sur $\smxs$, dont l'image est 
exactement $(I^l)_l^G$, partie homogène de degré $l$ de $(I^l)^G$.
\item{\rm (ii)} $R^q\pi_*(\es^l\V_L)|_{\es^m_*X}=0$ pour $q>0$.
\end{itemize}
\end{theor}

\par{\bf Preuve:}

Le foncteur image directe invariante $q_*^G$ est défini comme il suit
: pour un \fs $F$ sur $B$, et un ouvert $U$ de $\hilxs$,
$q_*^G(F)(U)=(F(q^{-1}(U)))^G$, ce qui a un sens puisque  $q^{-1}(U)$ est un ouvert $G$-invariant de $B$. En utilisant les propriétés des variétés quotient par un groupe 
fini, on obtient $\O^G_B=q^G_*(\O_B)=\O_{\hilxs}$ d'où 
$q^G_*(\es^l\tV)=\es^l\V_L$.

On en déduit que 
$\pi_*(\es^l\V_L)=\pi_*q_*^G(\es^l\tV)=p_*^G\rho_*(\es^l\tV)$ qui 
se plonge canoniquement dans $p_*^G(\es^lW)$, avec 
$p^G_*((I^l)_l)=(I^l)^G_l $ pour image. Les morphismes $p$ et $q$ sont 
finis donc leurs images directes supérieures sont nulles. Par 
composition des foncteurs dérivés $R\pi_*, Rp_*^G=p_*^G, R\rho_*$ 
et $Rq_*^G=q_*^G$ on trouve que
$R\pi_*(\es^l\V_L)=R\pi_*(q^G_*(\es^l\tV))=R\pi_*\circ 
Rq_*^G(\es^l\tV)=Rp_*^G\circ R\rho_*(\es^l\tV)=p_*^GR\rho_*(\es^l\tV)$ 
d'où la nullité des $R^q\pi_*(\es^l\V_L)$ pour $q>0$ sur 
$\smxs$.$\Box$

\begin{cor}
\begin{itemize}
\item{\rm (i)} 
$\H^0(\hil,\es^l\V_L\tens\de^{\otimes s})=\H^0(\ppm,(I^l)_l\tens\O(s,\ldots,s))^G$
\item{\rm (ii)} 
$\H^q(\hils,\es^l\V_L\tens\de^{\otimes s})=\H^q(\ppms,(I^l)_l\tens\O(s,\ldots,s))^G$
 pour $q>0$.
\end{itemize}
\end{cor}

\par {\bf Preuve:}

Compte-tenu du fait que $\de=\pi^*(\O(1,\ldots,1)^G)$ il suffit 
d'écrire 
$R^q\pi_*(\es^l\V_L\tens\de^{\otimes s})=R^q\pi_*(\es^l\V_L)\tens\O(s,\ldots,s)^G=((I^l)_l\tens\O(s,\ldots,s))^G$ si $q=0$ et $0$ sinon, et d'utiliser les 
propriétés de \coh locale pour (i) et la suite spectrale de Leray 
pour (ii).$\Box$

\begin{cor}
\label{h1v}
Le \th \ref{teore} est vrai pour $n\le 11$.
\end{cor}
\begin{rem}
{\rm La démonstration utilise le corollaire \ref{grinv} qui sera
  démontré au paragraphe \ref{1cap3h}, mais je préfère la donner ici pour motiver le travail fait dans le chapitre \ref{1cap4}}.
\end{rem}
\par {\bf Preuve:} Regardons le cas $l=1$. On rappelle que $E=\H^0(\pp,\O(1))$. Pour $L=\O_{\pp
}(2l-3)$ il faut  
calculer $\H^0(\ppm,V\tens\O(1,\ldots,1))^G$ où $V$ est le noyau du 
\mor \sur $W=L_1\oplus\cdots\oplus 
L_m\surto\L=\oplus_{i<j}L_{\Delta_{ij}}$. On tensorise la suite exacte 
$0\to V\to W\to\L\to 0$ par $\O(1,\ldots,1)$ et on écrit la suite 
exacte de \coh équivariante. On a vu que $\L$ n'avait pas de \coh 
$\sigm$-équivariante. On obtient ainsi 
$$\H^0(\ppm ,V\tens\O(1,\ldots,1))^G=\H^0(\ppm,\oplus_i\O(1,\ldots,2l-2,1,\ldots
,1))^G=$$
\begin{equation}
=\H^0(\pp,\O(2l-2))\tens\es^{m-1}\H^0(\pp,\O(1))=\es^{2l-2}E\tens\es^{m-1}E.
\label{AAA}
\end{equation}

On trouve aussi 
$\H^q(\ppm,V\tens\O(1,\ldots,1))^G=\H^q(\ppms,V\tens\O(1,\ldots,1))^G$ 
pour $q\le 2$ (puisque le complémentaire de l'ouvert $\ppms$ dans 
$\ppm$ est de codimension $4$) et comme le membre de gauche est nul pour 
$q\ge1$, celui de droite est nul pour $q=1$ et $q=2$. Donc 
$\H^1(\hils,\V_k\tens\de)=0$ d'où aussi $\H^1(\hil,\V_k\tens\de)=0$. 
À partir de la présentation (\ref{*}) de $\R$, avec $k$ remplacé par 
$2l-3=2\cdot 1-3=-1$, et $m=n+l^2=n+1$, on obtient 
$\V_{\O(-1)}\simeq\R$ et 
$$\dim H^0(\R\tens\de)=\dim\H^0(\V_{\O(-1)}\tens\de)=\dim\es^nE=\frac{(n+1)(n+2
)}{2}$$
pour $n$ tel que $3\le n\le 5$.

Passons ensuite à $l=2$. Il faut calculer cette fois-ci 
$\H^0(\ppm,(I^2)_2\tens\O(1,\ldots,1))^G$. On regarde les suites exactes 
associées à la \fil de $\es^2W$:
$$0\to I_2\to \es^2W\to\gr_0(\es^2W)\to 0$$
$$0\to (I^2)_2\to I_2\to\gr_1(\es^2W)\to 0$$

On sait que $\gr_0(\es^2W)=\es^2\L$. On  démontrera plus tard 
(cor.\ref{grinv}) que  $\gr_i(\es^lW)$ n'a pas de \coh équivariante 
si $l-i$ est impair.

En écrivant les suites exactes de \coh équivariante, après avoir 
tensorisé par $\O(1,\ldots,1)$, on obtient
$$\H^q((I^2)_2\tens\O(1,\ldots,1))^G\simeq 
\H^q(I_2\tens\O(1,\ldots,1))^G,\forall q\ge 0,$$
et que $\H^0(I_2\tens\O(1,\ldots,1))^G$ est le noyau du \mor
$$mor: \H^0(\es^2W\tens\O(1,\ldots,1))^G\to \H^0(\es^2\L\tens\O(1,\ldots,1))^G.$$
À l'aide de la proposition \ref{propinv} ces espaces d'invariants se calculent aisément pour donner
$$\H^0(\es^2W\tens\O(1,\ldots,1))^G=\es^{2(2l-3)+1}E\tens\es^{m-1}E\oplus\es^2(\es^{2l-3+1}E)\tens\es^{m-2}E$$
$$\H^0(\es^2\L\tens\O(1,\ldots,1))^G=\es^{2(2l-3)+2}E\tens\es^{m-2}E.$$
 
La composante de $mor$ sur le second facteur est induite par le \mor de restriction 
$$\H^0(\pp\times\pp,\O(2l-3+1)\boxtimes\O(2l-3+1))\to\H^0(\O_D(2(2l-3)+2)).$$
Le \mor $mor$ est alors \sur et on peut calculer la dimension de son noyau. L'espace
$\H^1((I^2)_2\tens\O(1,\ldots,1))^G$ s'injecte dans 
 $\H^1(\es^2W\tens\O(1,\cdots,1))^G=0$: il est donc 
nul. On obtient donc aussi la nullité de $\H^1(\hil,\es^2\V_L\tens\de)$.

À partir de la
présentation (\ref{**}) de $\es^2\R\tens\de$, et des annulations de la \coh
supérieure obtenues, il résulte une suite exacte de \reps
$$0\to\Lambda^2E\tens\H^0(\de)\to E\tens
\H^0(\V_1\tens\de)\to\H^0(\es^2\V_1\tens\de)\to\H^0(\U,\es^2\R\tens\de)\to
0$$
On en déduit
$$\dim\H^0(U,\es^2\R\tens\de)=\dim\es^3E\tens\es^{n+3}E\oplus\es^2(\es^2E)\tens\es^{n+2}E$$
$$-\dim\es^4E\tens\es^{n+2}E-\dim E\tens\es^2E\tens\es^{n+3}E+\dim\Lambda^2E\tens\es^{n+4}E.$$

Par suite,
$$\dim\H^0(U,\es^2\R\tens\de)=10\left(\!\begin{array}{c}n+5\\ 2\end{array}\!\right)+21\left(\!\begin{array}{c}n+4\\ 2\end{array}\!\right)$$
$$-15\left(\!\begin{array}{c}n+4\\ 2\end{array}\!\right)-
18\left(\!\begin{array}{c}n+5\\ 2\end{array}\!\right)+3\left(\!\begin{array}{c}n+6\\ 2\end{array}\!\right)=\frac{1}{2}(n+1)(n+2).$$
On a ainsi démontré le \th \ref{teore} pour tout $n$ tel que $3\le n\le 11$.$\Box$

Les vraies difficultés apparaissent à partir de $l=3$. On écrit à
nouveau les suites exactes associées à la \fil de $\es^3W$:
$$0\to I_3\to\es^3W\to\es^3\L\to 0 $$
$$0\to (I^2)_3\to I_3\to\gr_1(\es^3W)\to 0 $$
$$0\to (I^3)_3\to (I^2)_3\to\gr_2(\es^3W)\to 0 $$

Comme on le verra dans le corollaire \ref{grinv}, $\gr_2(\es^3W)$ n'a pas de \coh équivariante. On le
savait déjà  pour $\es^3\L$. Les suites de \coh équivariante 
associées nous fournissent alors:
$$\H^q((I^3)_3\tens\O(1,\ldots,1))^G\simeq\H^q((I^2)_3\tens\O(1,\ldots,1))
^G$$ 
pour $q\ge 0$ et
$$\H^q(I^3\tens\O(1,\ldots,1))^G\simeq\H^q(S^3W\tens\O(1,\ldots,1))^G$$ 
pour $q\ge 0$.

Alors l'espace recherché $\H^0((I^3)_3\tens\O(1,\ldots,1))^G$ 
s'obtient comme le noyau du \mor
$$\alpha:\H^0(S^3W\tens\O(1,\ldots,1))^G\to\H^0(\gr_1(S^3W)\tens\O(1,\ldots,1))^G$$
et puisque
$\H^1(S^3W\tens\O(1,\ldots,1))^G=0$, l'espace
$\H^1((I^3)_3\tens\O(1,\ldots,1))^G$ s'obtient comme son conoyau.

Toute la suite sera consacrée à l'étude minutieuse du \mor 
$\alpha$, afin de déterminer son noyau.

Pour bien comprendre la situation, on examinera d'abord le cas $m=2$, 
qui est essentiel pour pouvoir comprendre le cas $m$ général, dans 
un premier temps sans tensoriser avec le \fib inversible 
$\O(1,\ldots,1).$


\section{Le noyau du \mor $\alpha$}
\label{1cap4}

Dans les sections \ref{sec:leg}-\ref{sec:int} suivantes  on peut
 choisir $X$ comme {\'e}tant   une
surface lisse quasi-projective.

\subsection{Le gradu{\'e} $\gr_1(\es^3W)$}
\label{sec:leg}

On se propose de d{\'e}crire le gradu{\'e} $\gr_1(\es^3W)$, dans le cas 
$m=2$, $L$ \fib inversible sur $X$.  Dans ce cas les ouverts index{\'e}s par un {\'e}toile co{\"\i}ncident avec les 
espaces entiers, $D$ est la diagonale $\Delta$ de $X\times X$, 
$\L=L_D=L_{\Delta}$ et il n'y a pas de confusion si on le note 
toujours $L$. Ici $W=L_1\oplus L_2$.  Le gradu{\'e} $\gr_i(\sym W)$ est 
un $\sym L$-module; pour $i=0$, c'est l'alg{\`e}bre $\sym L$. On a vu 
aussi que $\gr_1(\sym W)$ est $\N_D$ en degr{\'e} $0$ et $V_D$ en degr{\'e} 
$1$. 

Pour comprendre $\gr_1(\es^kW)$ on regarde les $k+1$ \mors 
canoniques de $\es^kW|_D \to L^{\tens k}$ qui sont construits de la 
mani{\`e}re suivante: on consid{\`e}re les deux \mors canoniques 
$W|_D=L\oplus L\to L$ dont l'un, $\epsilon_+$, est donn{\'e} par la 
matrice $(id,id)$ et l'autre, $\epsilon_-$, est 
donn{\'e} par la matrice $(id,-id)$; c'est le \mor $\epsilon$ consid{\'e}r{\'e} au paragraphe \ref{filtration}. On obtient un isomorphisme 
$\varepsilon:W|_D\to L\oplus L$ d{\'e}fini par $(\epsilon_+,\epsilon_-)$ 
et par suite un isomorphisme
$$\es^k\varepsilon:\es^kW|_D\to\es^k(L\oplus L)=L^{\tens 
k}\oplus\cdots\oplus L^{\tens k}$$ 
dont la $i$-{\`e}me composante dans la 
somme directe est not{\'e}e $\varepsilon_{i,k-i}$. La derni{\`e}re 
composante $\varepsilon_{k,0}$ envoie $e_1^ie_2^{k-i}$ en 
$(-1)^{k-i}e^k$ avec pour $e,e_1=p_1^*(e),e_2=p_2^*(e)$  des rep{\`e}res locaux de 
$L,L_1$ et respectivement $L_2$, et s'{\'e}tend donc en un morphisme 
d'alg{\`e}bres $\sym W\to\sym L$ qui n'est autre que le \mor d'alg{\`e}bres 
consid{\'e}r{\'e} auparavant, $\epsilon_{-}$.

L'avant-derni{\`e}re composante $\varepsilon_{k-1,1}$ d{\'e}finit une 
d{\'e}rivation $\sym W\to\sym L$ compatible avec la graduation. En degr{\'e} 
$k$, $\varepsilon_{k-1,1}$ envoie $ e_1^ie_2^{k-i}$ sur 
$(i(-1)^{k-i}+(k-i)(-1)^{k-i-1})e^k$ et $\sym L$ est vu comme $\sym 
W$-module par l'interm{\'e}diaire du \mor $\epsilon_-$. On v{\'e}rifie alors que 
$$\varepsilon_{k-1,1}(xy)=\varepsilon_{k-1,1}(x)\epsilon_-(y)+\epsilon_-
(x)\varepsilon_{k-1,1}(y)$$ 
et comme $\epsilon_-(I)=0$ on obtient que le 
noyau de $\varepsilon_{k-1,1}$ contient $\F^2$. Donc 
$\varepsilon_{k-1,1}$ passe au quotient en une d{\'e}rivation lin{\'e}aire 
sur l'alg{\`e}bre $\sym L$, not{\'e}e encore $\epsilon_+:\gr_1(\sym 
W))\to\sym L$ qui est elle aussi compatible avec la graduation; cette 
propri{\'e}t{\'e}, jointe au fait qu'on conna{\^\i}t d{\'e}j{\`a} $\epsilon_+$ 
sur $V_D=\gr_1(W)$, caract{\'e}rise la d{\'e}rivation $\epsilon_+$.

Le \fs conormal {\`a} $D$, $\N_D$, est isomorphe au \fs $\Omega^1$ des formes diff{\'e}rentielles sur $X$. 
Un tel isomorphisme s'obtient en associant {\`a} la diff{\'e}rentielle ${\rm 
d}f$ d'une fonction r{\'e}guli{\`e}re sur un ouvert $U$ de $X$, la section 
de $\N_D$ d{\'e}finie par la classe $[f_2-f_1]$ o{\`u} $f_i=pr_i^*(f)$. 
L'image directe par $p$ de la suite (\ref{s3}), {\'e}crite en degr{\'e} $k$ 
est (c'est un cas particulier de la proposition \ref{3.1})~:
$$0\to \N\tens L^{\tens k}\to\gr_1(\es^kW)\to K\tens L^{\tens(k-1)}\to 0$$
et comme ici $K\simeq L$ on obtient une suite exacte de $\O_X$-modules
\begin{equation}
0\to\Omega^1\tens L^{\tens k}\to\gr_1(\es^kW)\to L^{\tens k}\to 0 
\label{s7}
\end{equation}
o{\`u} la premi{\`e}re fl{\`e}che se calcule de la mani{\`e}re suivante: pour 
$f\in\O(U)$ et $s\in\H^0(U,L^{\tens k})$, 
$${\rm d}f\tens s 
\mapsto\frac{1}{2}[(f_2-f_1)(s_1+(-1)^ks_2)]=(-1)^k[(f_2-f_1)s_2]$$
o{\`u} $s_i$ est la section de $\es^kW$ sur $U\times U$ d{\'e}finie par 
$pr_i^*(s)$. La seconde fl{\`e}che est $\epsilon_+$. En effet 
$\frac{1}{2}(s_1+(-1)^ks_2)$ est une section de $\es^kW$ dont l'image par 
$\epsilon_-$ est $s$ (et c'est aussi le cas pour $(-1)^ks_2$).


\subsection{Les opérateurs $\nabla$ et $\Delta$}

On considère l'opérateur $\comp$-linéaire $\nabla:L^{\tens 
k}\to\gr_1\es^kW$ qui associe à une section $s$ la classe $\nabla(s)$ 
de la section de $\F^1\es^kW$ définie par $(-1)^ks_2-s_1$; autrement dit

$\nabla(s)=[(-1)^ks_2-s_1]$
(on vérifie que $\epsilon_-(\nabla(s))=0$).

On a $\nabla(1)=[pr_2^*1-pr_1^*1]=0$.

\begin{lemme}

Cet opérateur n'est pas linéaire, mais satisfait à la condition 
$\nabla(fs)={\rm d}f\tens s+f\nabla(s)$ où $f$ est une fonction 
régulière sur un ouvert $U\subset X$ et $s$ est une section locale 
de $L^{\tens k}$ sur $U$. Dans cette formule $\Omega^1\tens L^{\tens 
k}$ est vu comme sous-module de $\gr_1\es^kW$ par l'inclusion de la suite 
(\ref{s7}).

\end{lemme}

\par{\bf Preuve:}

Par définition :
\begin{eqnarray*}
\nabla(fs)&=&[(-1)^kf_2s_2-f_1s_1]
=[(-1)^k(f_2-f_1)s_2+f_1((-1)^ks_2-s_1)]= \\ 
&=&{\rm d}f\tens s +f\nabla(s).
\end{eqnarray*}

Ceci signifie que $\nabla$ est un opérateur différentiel linéaire 
d'ordre $1$. Ce qu'on va utiliser c'est que $-\frac{\nabla}{2k}$ est une 
section $\comp$-linéaire de $\epsilon_+$ dans la suite (\ref{s7}), en vérifiant par un 
calcul direct que $-\frac{1}{2k}\epsilon_+(\nabla(s))=s$.

Au passage on peut remarquer que $\nabla$ se factorise en un \mor 
$\O_X$-linéaire, \hbox{$\bar{\nabla}:J^1L^k\to\gr_1(\es^kW)$} ($J^1L^k$ est le 
\fib des jets à valeurs dans $L^{\tens k}$, avec sa structure naturelle de 
$\O_X$-module à gauche) qui rend commutatif le diagramme
$$\begin{array}{ccccccccc}
{0}&\rightto{}&{\Omega_X^1\tens L^{\tens k}}&\rightto{}&{J^1L^k}&\rightto{}&{
  L^{\tens k}}&\rightto{}&{ 0} \\
&&\downto{id}&&\downto{{\bar{\nabla}}}&\swarrow_{\nabla}&\downto{-2k}& & \\
{0}&\rightto{}&{\Omega_X^1\tens L^{\tens k}}&\rightto{}&{\gr_1\es^kW}&\rightto{}&{L^{\tens 
k}}&\rightto{}&0
\end{array}$$

Les deux flèches extrêmes sont des isomorphismes, donc la flèche 
du milieu est un isomorphisme. Cet \isom induit en particulier, pour $k=1$, 
un \isom $J^1L\to V_D$.

\vs Considérons maintenant deux indices $i$ et $j$ tels que $i+j=k$. 
Soient $s\in\H^0(U,L^{\tens i})$ et $t\in\H^0(U,L^{\tens j})$ deux 
sections locales au-dessus du même ouvert $U$. On considère les 
sections $s_k$ de $\es^iW$ et $t_k$ de $\es^jW$ définies par image 
réciproque par $pr_k$ (pour $k=1,2$) et la section de $\gr_1\es^kW$:
$$D(s,t)=[(-1)^is_2t_1-(-1)^js_1t_2].$$

Comme $\epsilon_-(D(s,t))=0$, la quantité entre crochets appartient 
à $\F^1\es^kW$ et par conséquent la formule a un sens.

\begin{lemme}

Soient $s$ et $t$ comme ci-dessus. On a dans $\gr_1(\sym W)$, considéré comme \hbox{$\sym L$-module} :
$$\begin{array}{ccc}
\nabla(st)&=&\nabla(s)t+s\nabla(t)\\
D(s,t)&=&\nabla(s)t-s\nabla(t)\\
\end{array}$$
\end{lemme}

\par{\bf Preuve:}

Compte-tenu de la définition de l'homomorphisme $\epsilon_-:W|_D\to 
L$, on a dans $\sym W$:
\begin{eqnarray*}(-1)^k\nabla(st)&=&[s_2t_2-(-1)^ks_1t_1]=[(s_2-(-1)^is_1)t_2+(-1)^is_1(t_2-(-1
)^jt_1)] \\
&&=(-1)^k(\nabla(s)t+s\nabla(t))
\end{eqnarray*}

De même
$$(-1)^i[s_2t_1-(-1)^ks_1t_2]=[(-1)^i(s_2-(-1)^is_1)t_1]-(-1)^j[s_1(t_2-
(-1)^jt_1)]$$
compte-tenu de la définition de la structure multiplicative dans 
l'algèbre bigraduée $\gr(\sym W)$ et que $gr_0(\sym W)=\sym L$, ceci n'est autre que $\nabla(s)t-s\nabla(t)$.

\begin{cor}

Soit $k=i+j$, et $l=i-j$. L'opérateur différentiel $(s,t)\mapsto 
-l\nabla(st)+kD(s,t)$ prend ses valeurs dans $\Omega^1\tens L^{\tens 
k}$.

\end{cor}

Ceci résulte du fait que $\epsilon_+(\nabla(s))=-2ks$ si $s$ est une section locale de $L^{\tens k}$: parce que $\epsilon_+$ est une 
dérivation, ceci entraîne en effet que 
$\epsilon_+(-l\nabla(st)+kD(s,t))=0$.


\subsection{Le \mor $\alpha_2 :\H^0(X\times 
X,\es^kW)^{\tau}\to\H^0(X,\gr_1(\es^kW))^{\tau}$}

On suppose $k$ impair.
Soit $\tau$ la transposition $(12)$ et désignons par 
$\H^0(\es^kW)^{\tau}$ l'espace des sections de $\H^0(\es^kW)$ invariantes 
sous l'action de $\tau$. Puisque dans le cas où $k$ est impair, le 
gradué $\gr_0(\es^kW)$ n'a pas de \coh équivariante, ces sections 
invariantes définissent des sections de $\F^1\es^kW$, d'où le \mor  
$\alpha_2$. En outre, les sections de $\gr_1(\es^kW)$ sont invariantes pour
l'action de $\tau$, puisque pour $k-1$ pair, $\tau$ agit 
trivialement sur tous les gradués de sa filtration déduite de
(\ref{filtratio}) (voir cor. \ref{grinv}). On 
a un isomorphisme canonique 
$$\oplus_{i>j,i+j=k}\H^0(L^{\tens i})\tens_{\comp}\H^0(L^{\tens j})\simeq\H^0(\es^kW)^{\tau}$$
donné par $s\tens t\mapsto s_1t_2+s_2t_1$ pour $s\in\H^0(L^{\tens i})$ et $t\in\H^0(L^{\tens j})$.

Désignons par $\mu:\H^0(L^{\tens i})\tens \H^0(L^{\tens 
j})\to\H^0(L^{\tens k})$ la multiplication et considérons pour $k>1$ 
le scindage de la suite (\ref{s7}) sur les sections globales :

$$\H^0(\gr_1\es^kW)\simeq\H^0(\Omega^1\tens L^{\tens 
k})\oplus\H^0(L^{\tens k})$$
défini sur le premier facteur par l'inclusion canonique, et sur le 
deuxième facteur par la section \hbox{$s\mapsto\nabla(s)$.}

Si $s\in\H^0(X,L^{\tens i})\tens\H^0(X,L^{\tens j})$ on pose 
$\Delta_l(s)=D(s)-\frac{l}{k}\nabla\mu(s)$ pour $l=i-j$ et 
$i+j=k$. Pour $l=k$ on obtient que $\Delta_l(s)=0$ puisque pour 
$s$ décomposable en $m\tens 1$, 
$\Delta_l(s)=\nabla(m)\cdot1-m\cdot\nabla(1)-\nabla(m)=0$.

\begin{prop}

\label{propA}

La matrice de $\alpha_2$ dans ces  décompositions est donnée par

$$\left(\begin{array}{ccccc}
0&\cdots &(-1)^{\frac{k+l}{2}}\Delta_l & \cdots&(-1)^{\frac{k+1}{2}}\Delta
_1 \\
2k &\cdots &-(-1)^{\frac{k+l}{2}}2l\mu
&\cdots&-(-1)^{\frac{k+1}{2}}2\mu
\end{array}\right).$$

\end{prop}

\par{\bf Preuve :}

On part d'une section $S\in\H^0(\es^kW)^{\tau}$ qui provient d'une 
section de $\H^0(X,L^{\tens i})\tens\H^0(X,L^{\tens j})$ avec $i+j=k$, 
$i-j=l$. Supposons aussi que $j\ne 0$ et que cette section  se 
décompose en $s\tens t$ avec $s\in\H^0(L^{\tens i})$ et 
$t\in\H^0(L^{\tens j})$. Alors $S$ s'écrit dans $\H^0(\es^kW)^{\tau}$ 
comme $s_1t_2+t_1s_2$ et son image par $\alpha_2$ dans $\H^0(\gr_1\es^kW)$ est 
la classe $[ s_1t_2+t_1s_2]$ modulo $\F^2\es^kW$.

Si $i$ est pair et $j$ impair, $[ s_1t_2+t_1s_2]=D(s,t)$ et la 
composante dans $\H^0(L^{\tens k})$ est $\epsilon_+(D(s,t)).$ Mais 
$\epsilon_+(-l\nabla(st)+kD(s,t))=0$ donc 
$\epsilon_+(D(s,t))=\frac{l}{k}\epsilon_+(\nabla(st))=\frac{l}{k}\cdot(-2k\mu(s\tens t))=-2l\mu(s\tens t)$. Si $i$ est impair on obtient 
l'opposé. Si $l=k$ alors $j=0$, $s\in\H^0(L^{\tens k})$ et son 
écriture dans $\H^0(\es^kW)^{\tau}$ est $s_1+s_2$. Son image par $\alpha_2$ 
est $[s_1+s_2]=-\nabla(s)$ et $\epsilon_+(-\nabla(s))=2ks$.  D'où 
la deuxième ligne de la matrice.

Pour trouver la composante dans $\H^0(\Omega^1\tens L^{\tens k})$ il 
suffit de soustraire l'image réciproque par la section 
$-\frac{1}{2k}\nabla$ du \mor $\epsilon_+$ de la composante dans 
$\H^0(L^{\tens k})$. Par exemple pour $i$ pair on fait 
$D(s,t)-(-\frac{1}{2k}\nabla(\epsilon_+(D(s,t))))=D(s,t)+\frac{1}{2k}\nabla(-2ls
t)=\Delta_l(S)$. Pour $i$ 
impair on trouve l'opposé et pour 
$j=0$ : $-\nabla(s)+\frac{1}{2k}\nabla(2ks)=0$, d'où la première 
ligne.$\Box$

\begin{cor}

Si $k=3$, le noyau et le conoyau de $\alpha_2$ sont isomorphes respectivement 
au noyau et au conoyau de l'opérateur linéaire

$$3D-\nabla\mu:\H^0(L^{\tens 2})\tens \H^0(L)\to\H^0(\Omega^1\tens 
L^{\tens 3}).$$

\end{cor}

\par{\bf Preuve :}

La matrice de $\alpha_2$ s'écrit ici

$$\left(\begin{array}{cc}
0&\Delta_1=D-\frac{1}{3}\nabla\mu\\
6&-2\mu
\end{array}\right)$$
ce qui conduit immédiatement à l'énoncé.$\Box$


\subsection{Généralisation à $\hilx$}
\label{1cap3h}

Le cas général repose essentiellement sur le cas $m=2$. Il faut 
considérer les invariants par rapport au groupe symétrique 
$G=\sigm$ mais on a donné dans les préliminaires, paragraphe \ref{calinv}, le procédé qui nous 
ramène à des calculs d'invariants plus aisés.

\subsubsection{{Description de } $\H^0(\xms,\gr_1(\es^kW))^G$}

Le \fs $\gr_1(\sym W)$ a pour support la diagonale $D$ de $\xms$. Le 
groupe symétrique agit sur la situation. Soit $U_{12}$ le 
complémentaire de la réunion des diagonales $\Delta_{i,j}$ pour $\{i,j\}\ne\{1,2\}$ dans 
$\xms$. Cet ouvert contient uniquement la diagonale 
$\Delta_{1,2}=\Delta\times X^{m-2}\cap\xms$.  On note 
$W_{12}=(L_1\oplus L_2)|_{U_{12}}$ et 
$W^{12}=(L_3\oplus\cdots\oplus L_m )|_{U_{12}}$ de sorte que $W 
|_{U_{12}}=(W_{12}\oplus W^{12})|_{U_{12}}$.

\begin{prop}

\label{restr}

On a 
$$\gr_k(\sym W)|_{\Delta_{12}}\simeq\oplus_{i+j=k}\gr_i(\sym 
W_{12})\tens\es^jW^{12}[-j].$$
\end{prop}

\par{\bf Preuve :}

Le produit tensoriel peut être vu comme un produit tensoriel sur 
$\comp$, ou bien comme un produit tensoriel externe, si, plutôt 
qu'utiliser $W_{12}$ et $W^{12}$ on utilise $W'_{12}=L_1\oplus L_2$ 
sur $X\times X$ et\ \ \ 
 $W^{12\prime}=L_3\oplus\cdots\oplus L_m$ sur $X^{m-2}$. 
On a  $W_{12}=pr_{12}^*(W'_{12})|_{U_{12}}$, 
$W^{12}=pr_{3\cdots m}^*(W^{12\prime})|_{U_{12}}$ et $W'_{12}\boxtimes W^{12\prime}=W_{12}\tens W^{12}$. La notation 
$\es^jW^{12}[-j]$ signifie qu'on place $\es^jW^{12}$ en degré $j$.

Le \mor $\varepsilon$ devient en restriction à $U_{12}$ :
$$\epsilon|_{U_{12}}:\sym W |_{U_{12}}\to\sym L|_{U_{12}}=\sym 
L_{\Delta}|_{\Delta_{12}}$$
et, puisque $\sym W |_{U_{12}}=\sym W_{12}\tens\sym W^{12}=\sym W'_{12}\boxtimes\sym W^{12\prime}$, ce \mor 
est aussi un produit tensoriel externe des \mors
$$\epsilon_{12}:\sym W'_{12}\to\sym L_{\Delta}$$
de noyau $I_{12}$, qui recopie la situation étudiée dans le cas où 
$m$ était égal à $2$, et
$$\epsilon_{3\cdots m}:\sym W^{12\prime}\to\O_{X^{m-2}}$$
qui vaut l'identité en degré $0$ et $0$ en degré $\ge 1$, de noyau 
$(\sym W^{12\prime})_{\ge 1}$.

On filtre $\sym W'_{12}$ par les puissances de $I_{12}$ et $\sym
W^{12\prime}$ par les puissances de $(\sym W^{12\prime})_{\ge 1}$ qui
sont égales aux $(\sym W^{12\prime})_{\ge j}$.

Le noyau de $\epsilon|_{U_{12}}$ s'écrit alors comme
$$I=I_{12}\boxtimes\sym W^{12\prime}+\sym W'_{12}\boxtimes\sym W^{12\prime}_{\ge 1}$$
et sa puissance $k$-ième
$$I^k=\sum_{i+j=k}I^i_{12}\boxtimes \sym W^{12\prime}_{\ge j}.$$

On veut calculer $I^k/{I^{k+1}}.$ Le calcul de ce gradué se fait à
l'aide du lemme \ref{altalema} du paragraphe préliminaire. La condition cohomologique d'annulation est
vérifiée en vertu du lemme préliminaire \ref{incaolema}.

On trouve 
$$\gr_k(\sym W)|_{\Delta_{12}}\simeq 
I^k/{I^{k+1}}|_{U_{12}}\simeq\oplus_{i+j=k}\gr_i(\sym 
W_{12})\tens\es^jW^{12}[-j].\Box$$

\begin{cor}
On a un \iso
$$\H^0(\xms,\gr_k(\sym W))^G\simeq\oplus_{i+j=k}(\H^0(X^2,\gr_i(\sym 
W_{12}))^{\sigma_2}\tens\H^0(X^{m-2},\es^jW^{12})^{\sigma_{m-2}}[-j]).$$

\end{cor}

\par{\bf Preuve du corollaire:}

Appliquons le résultat du lemme \ref{propinv} pour 
$M=\H^0(\xms,\gr_k(\sym W))$ et pour l'ensemble d'indices $I=\{\{i,j\}\}_{1\le i<j\le m}$ sur 
lequel $G$ agit. Prenons $L_{1,2}=\gr_k(\sym W)|_{U_{12}}$, calculé 
par la proposition \ref{restr},  et $L_{i,j}$ le \fib similaire sur la 
diagonale $\Delta_{ij}$: $\gr_k(\sym W)|_{\Delta_{ij}}$. L'espace $M^G$ 
s'obtient en prenant les invariants de $M_{{1,2}}$, pour le 
stabilisateur de $\{1,2\}$,
$\Stab\{1,2\}=\sigma_2\times\sigma_{m-2}$. Le complémentaire de
l'ouvert $\Delta_{i,j}$ dans $\Delta\times X^{m-2}$ est de codimension
$\ge 2$, donc pour le calcul de l'espace des sections $\H^0(\gr_i(\sym W_{12})\tens \es^j 
W^{12}[-j])$ on peut se placer sur $\Delta\times X^{m-2}$, où on
applique le \th de Künneth.
Puisque  $\sigma_2$ n'agit pas sur ce qui provient de $X^{m-2}$ 
et $\sigma_{m-2}$ n'agit pas sur ce qui provient de $X^2$, on a
$$\H^0(\gr_i(\sym W_{12})\tens \es^j 
W^{12}[-j])^{\sigma_2\times\sigma_{m-2}}=\H^0(\gr_i(\sym 
W_{12}))^{\sigma_2}\tens\H^0(\es^jW^{12})^{\sigma_{m-2}}[-j]$$
d'où le résultat.

\begin{cor}

\label{grinv}
Pour $l-k$ impair, le \fs $\gr_k(\es^lW)$ n'a pas de sections invariantes 
sous l'action de $G$.
\end{cor}

\par{\bf Preuve:}

En effet $\gr_i(\es^lW_{12})$ n'a pas de sections invariantes sous l'action de 
$\sigma_2$ si $l-i$ est impair puisque dans la \fil déduite de
(\ref{filtratio}), aucun de ses gradués n'a des sections invariantes
($\sigma_2$ agit par $(-1)$ sur le \fib conormal de la diagonale dans
$X\times X$ et sur $L$, et trivialement sur $K$ donc par $(-1)^{2q+l-i}$ sur $\gr_q(\gr_i(\es^lW_{12}))$). Mais
$$\gr_k(\es^lW)\simeq\oplus_{i+j=k}\gr_i(\es^{l-j}W_{12})\tens\es^jW^{12}$$
et donc si $l-i-j=l-k$ est impair, $\gr_i(\es^{l-j}W_{12})$ n'a pas de 
sections invariantes sous l'action de $G$. On avait déjà utilisé 
ce corollaire dans la section \ref{calco}.$\Box$

\vs On a tout fait pour comprendre que pour $l$ impair
\begin{eqnarray}
\label{nume1}
\H^0(\gr_1\es^lW)^G&=&\H^0(\gr_1\es^lW_{12})^{\sigma_2}\oplus\left[\H^0(\gr_0\es
^{l-1}W_{12})^{\sigma_2}\tens\H^0(W^{12})^{\sigma_{m-2}}\right]\nonumber\\
&=&\H^0(J^1L^{\tens l})\oplus\left[\H^0(L^{\tens(l-1)})\tens\H^0(L)\right]\nonumber\\
&=&\H^0(\Omega^1\tens L^{\tens k})\oplus\H^0(L^{\tens k})\oplus\H^0(L^{\tens(l-1)})\tens\H^0(L).
\end{eqnarray}

Les invariants de $\H^0(\xm,\es^lW)$ se calculent facilement mais 
l'écriture est lourde pour $l$ élevé. On préfère donc se 
limiter dans la suite au seul cas qui nous intéresse $l=3$.

On a 
\begin{equation*}
\label{}
\es^3W=\oplus^n_{i=1}L_i^{\tens 3}\oplus\oplus_{i\ne j}(L_i^{\tens 
2}\tens L_j)\oplus\oplus_{i<j<k}(L_i\tens L_j\tens L_k)
\end{equation*}
et la même proposition appliquée à $I=\{1,2,\ldots,m\}$ et les 
\fibs $L_i^{\tens 3}$, ensuite à $I=\{(i,j)\}_{1\le i,j\le m}$ et 
$L_{(i,j)}=L_i^{\tens 2}\tens L_j$ et finalement à $I=\{\{i,j,k\} 
\}_{1\le i<j<k\le m}$ et $L_{\{i,j,k\}}=L_i\tens L_j\tens L_k$, nous prouve 
que
\begin{equation}
\label{nume}
\H^0(X,L^{\tens 3})\oplus\left[\H^0(X,L^{\tens 
2})\tens\H^0(X,L)\right]\oplus\es^3\H^0(X,L)=\H^0(\xm,\es^3W)^{\sigm}
\end{equation}
l'\iso étant donné par
$$(s,t\tens u,vzw)\mapsto(\sum_i s_i,\sum_{1\le i<j\le m}(t_iu_j+t_ju_i),\sum_{i\ne j\ne k,i\ne k}v_iz_jw_k)$$
(naturellement $s_i=pr_i^*(s)$, et de même pour $t_i$, $u_i$, $v_i$, $z_i$, $w_i$).

En effet, comme $\H^0(\xm,L_i^{\tens 3})=\H^0(X,L^{\tens 3})$ par la 
formule de Künneth (car $L_i=pr_i^*L\tens(\tens_{j\ne i}pr^*_j\O_X)$) toutes les sections de $L_i^{\tens 3}$ sur $\xm$ sont en 
effet des images réciproques $pr_i^*(s)$ avec $s$ section de $L^{\tens 
3}$ sur $X$. Comme $\H^0(\xm,\oplus^m_{i=1}L_i^{\tens 3})^{\sigm}=\H^0(\xm,L_1)^{\sigma_{m-1}}$ et que le stabilisateur de 
$1$ n'agit pas sur $pr_1^*(s)$ on obtient que $\H^0(\xm,L_1^{\tens 
3})^{\sigma_{m-1}}=\H^0(X,L^{\tens 3})$ et de manière analogue les 
autres termes dans la décomposition. Le terme $\es^3\H^0(X,L)$ 
s'obtient puisqu'il faut considérer le stabilisateur de $\{1,2,3\}$ en 
tant qu'ensemble, c'est-à-dire $\sigma_3\times\sigma_{m-3}$ et

$(\H^0(\xm,L_1)\tens\H^0(\xm,L_2)\tens\H^0(\xm,L_3))^{\sigma_3}=\es^3\H^
0(X,L)$.

On dispose comme précédemment des opérateurs
$$\nabla:L^{\tens 
l}\to\gr_1(\es^lW)$$
défini par $\nabla(s)=\sum_{i<j}\nabla_{ij}(s)$ et
$$D:\Gamma(U,L^{\tens p})\times\Gamma(U,L^{\tens 
q})\to\Gamma(U\times\cdots\times U,gr_1(\es^lW))$$
défini pour $p+q=l$ et $U$ ouvert de $X$ par 
$D(s,t)=\sum_{i<j}D_{ij}(s,t)$.

Ici $\nabla_{ij}$ et $D_{ij}$ sont définis sur les ouverts $U_{ij}$ 
contenant la seule diagonale $\Delta_{ij}$ exactement comme $\nabla$ et 
$D$ dans le cas $m=2$, et jouissent des même propriétés:
$$\nabla_{ij}:L^{\tens 
l}\to\gr_1\es^lW$$
défini par $s\mapsto [-s_i+(-1)^ls_j]$ où $s_i=pr^*_i(s)$ et
$$D_{ij}:\Gamma(U,L^{\tens p})\times \Gamma(U,L^{\tens 
q})\to\Gamma((U\times\cdots\times U)\cap U_{ij},\gr_1\es^lW)$$
défini par $(s,t)\mapsto[(-1)^ps_it_j-(-1)^qs_jt_i]$ où $t_i=pr_i^*(t)$.

La proposition qui suit est l'analogue de la proposition \ref{propA}:

\begin{prop}
\label{defnu}
Dans les sommes directes (\ref{nume1}) et (\ref{nume}), la matrice du \mor canonique:
$$\alpha_m:\H^0(\xm,\es^3W)^G\to\H^0(\xm,\gr_1(\es^3W))^G$$
s'écrit
$$\left(\begin{array}{ccc}
0&\Delta&0\\
6&-2\mu&0\\
0&2id &-2\nu
\end{array}\right)$$
où $\nu$est le \mor canonique $\es^3\H^0(L)\to\H^0(L^{\tens 
2})\tens\H^0(L)$ induit par l'application linéaire $stu\mapsto st\tens 
u+su\tens t+ut\tens s$. Dans ce contexte $\Delta=D-\frac{1}{3}\nabla\mu$.
\end{prop}

\par{\bf Preuve:}
\vs Afin de calculer la première colonne de la matrice de $\alpha_m$, 
considérons une section locale $s$ de $L^{\tens 3}$ sur un ouvert 
$U$. La section définie par $s_1+s_2+\cdots+s_m$ est 
$\sigm$-invariante. C'est une section $\sigm$-invariante de $\F^1\es^3W$, 
ou bien une section $\sigma_2\times\sigma_{m-2}$-invariante de 
$\F^1\es^3W|_{U_{12}}=\F^1\es^3W_{12}\oplus\es^2W_{12}\tens W^{12}\oplus 
W_{12}\tens\es^2W^{12}\oplus\es^3W^{12}$. Modulo $\F^2\es^3W|_{U_{12}}$ on 
obtient $[s_1+s_2]\in\gr_1(\es^3W_{12})$ (puisque 
$s_3+\cdots+s_m\in\es^3W^{12}$ qui est inclus dans 
$\F^2\es^3W|_{U_{12}}$).

Son image dans la décomposition de $\H^0(\gr_1\es^lW))^G$ est 
$(-\nabla(s),0)$. Si on décompose encore
$\H^0(J^1L^3)=\H^0(\Omega^1\tens L^{\tens 3})\oplus\H^0(L^{\tens 3})$, 
d'après le résultat trouvé dans le cas $m=2$, on obtient la 
première colonne de la matrice comme $(0,6,0)$.

\vs Pour la deuxième colonne, considérons deux sections 
$s\in\Gamma(U,L^{\tens 2})$ et $t\in\Gamma(U,L)$. La section définie 
par $\sum_{1\le i<j\le 
m}(s_it_j+s_jt_i)=s_1t_2+s_2t_1+(s_1+s_2)(t_3+\cdots+t_m)+(t_1+t_2)(s_
3+\cdots+s_m)+\sum_{3\le i<j\le m}(s_it_j+s_jt_i)$ est 
$\sigm$-invariante. On procède comme auparavant. Modulo 
$\F^2\es^3W|_{U_{12}}$ , il reste seulement les deux premiers termes de 
cette expression :$D_{12}(s,t)=[s_1t_2+s_2t_1]\in\gr_1\es^3W_{12}$ et 
$[s_1+s_2]t\in\gr_1(\es^2W_{12})\tens (W^{12})^{\sigma_{m-2}}$. La classe 
$[s_1+s_2]$ modulo $\F^2\es^2W_{12}=\es^2I_{12}$ est son image dans 
$\es^2L_{\Delta}=L_{\Delta}^{\tens 2}$ soit $2s$. Au total, en 
utilisant aussi la décomposition de $D(s,t)$ trouvée dans le cas 
$m=2$ on obtient $(\Delta,-2\mu,2id)$.

Finalement, la troisième colonne s'obtient en partant de trois 
sections $s,t,u$ de $L$ sur un ouvert $U$. La section $\sum_{i\ne j\ne k, i\ne k}s_it_ju_k$ s'écrit comme
$$(s_1t_2+s_2t_1)(u_3+\cdots+u_m)+(s_1u_2+s_2u_1)(t_3+\cdots+t_m)+(t_1u_2+
t_2u_1)(s_3+\cdots+s_m)$$
$$+\left[(s_1+s_2)\sum_{3\le i\ne j\le 
m}t_iu_j\right]+\left[(t_1+t_2)\sum_{3\le i\ne j\le m}s_iu_j\right]$$
$$+\left[(u_1+u_2)\sum_{3\le i\ne j\le m}s_it_j\right]+\sum_{3\le i\ne l\ne 
k\le m,i\ne k}s_it_ju_k$$ 
et elle est $\sigm$-invariante. La classe dans 
$\gr^1\es^3W|_{U_{12}}$ de sa restriction à $U_{12}$ est 
$[(s_1t_2+s_2t_1)(u_3+\cdots+u_m)]+[(s_1u_2+s_2u_1)(t_3+\cdots+t_m)]+[(t_
1u_2+t_2u_1)(s_3+\cdots+s_m)]$ et chacune de ces composantes appartient 
à $\es^2W_{12}\tens W^{12}$. Le reste appartient à 
$\F^2\es^3W|_{U_{12}}$. Pour trouver leurs images dans 
$\gr_0(\es^2W_{12})\tens L$ on regarde les images de $s_1t_2+s_2t_1$, 
$s_1u_2+s_2u_1$ et $t_1u_2+t_2u_1$ dans $L_{\Delta}^{\tens 2}$ par le 
morphisme $\es^2W_{12}\to L^{\tens 2}_{\Delta}$. Or $s_it_j\mapsto -st$. La 
troisième colonne s'écrit $(0,0,-2\nu)$ où $\nu$ est le \mor 
canonique $\es^3\H^0(L)\to\H^0(L^{\tens 2})\tens \H^0(L)$ induit par 
l'application linéaire $stu\mapsto st\tens u+su\tens t+ut\tens s$. 
$\Box$

\begin{cor}

Pour $m\ge 3$, l'espace vectoriel des sections $\H^0(\es^3\V_L)$ sur 
l'ouvert $\hilxs$ est isomorphe à $\es^3\H^0(L)$ et l'espace vectoriel 
de \coh $\H^1(\es^3\V_L)$ est isomorphe à $\H^0(\Omega^1\tens L^{\tens 
3})$.

\end{cor}

\par{\bf Preuve du corollaire:}

Les espaces considérés sont le noyau et respectivement le conoyau du 
\mor $\alpha_m$. Il faut voir qu'ils coïncident avec le noyau et 
respectivement le conoyau du \mor nul:
$$\Delta\nu:\es^3\H^0(L)\to\H^0(\Omega^1\tens L^{\tens 3})$$

En effet, si $s$ est une section de $L$, on a $\nu(s^3)=3s^2\tens s$ 
et l'image de cette classe par $\Delta $ est 
nulle: $3D(s^2,s)-\nabla(s^3)=3\Delta(s^2)s-3s^2\nabla(s)-\nabla(s^2)s-s
^2\nabla(s)=2\nabla(s^2)s-4s^2\nabla(s)=0$. Comme ces sections 
engendrent $\es^3\H^0(L)$, ceci prouve que le \mor $\Delta\nu$ est nul.

On voit que, si $(a,b,c)\in \H^0(\xm,\es^3W)^G$, alors
$$\left(\begin{array}{ccc}
0 & \Delta&0\\
6&-2\mu&0\\
0&2id&-2\nu
\end{array}\right)
\left(\begin{array}{c}
a\\
b\\
c
\end{array}\right)
=\left(\begin{array}{c}
0\\
0\\
0
\end{array}\right)
\Longleftrightarrow
\begin{cases}
\Delta b=0\cr
6a-2\mu b=0\cr
2b-2\nu c=0\cr
\end{cases}$$
$$\Longleftrightarrow
\begin{cases}
a=\frac{\mu b}{3}\cr
\Delta\nu c=0\cr
b=\nu c\cr
\end{cases}$$

Par suite  $\Ker \alpha_m\simeq\Ker\Delta\nu\simeq\es^3\H^0(L)$ et

$$\left(
\begin{array}{ccc}
0&\Delta&0\\
6&-2\mu&0\\
0&2id&-2\nu\\
\end{array}
\right)
\left(
\begin{array}{c}
a\\
b\\
c\\
\end{array}
\right)
=
\left(
\begin{array}{c}
a'\\
b'\\
c'\\
\end{array}
\right)
\Longleftrightarrow
\begin{cases}
\Delta b=a'\cr
6a-2\mu b=b'\cr
2b-2\nu c=c'\cr
\end{cases}$$
 
$$\Longleftrightarrow
\begin{cases}
a=\frac{b'+2\mu b}{6}\cr
2a'-2\Delta\nu c=\Delta c'\cr
2b-2\nu c=c'\cr
\end{cases}$$

Par suite $(a',b',c')\in\Im\alpha_m\Longleftrightarrow 2a'-\Delta c'\in\Im\Delta\nu$ et donc $\coker\alpha_m\simeq\coker\Delta\nu$. Compte-tenu de la remarque préliminaire, $\coker\alpha_m\simeq\H^0(\Omega^1\tens L^{\tens 3})$.
$\Box$


\subsection{Introduction du \fib déterminant}
\label{sec:int}

On considère un autre \fib inversible $A$ sur $X$ auquel est associé 
un \fib vectoriel $\A=\boxtimes_iA_i$ sur $\xm$ et un \fib inversible 
quotient $\A/\sigm$ sur $\es^mX$. On désigne par $\D$ le \fib image 
réciproque de  $\A/\sigm$ sur $\hilx$ par le \mor de Hilbert-Chow $\hilx\to\es^mX$. Le 
problème est de déterminer l'espace vectoriel des sections de 
$\es^3V\tens\D$. Le calcul des invariants de $\H^0(\es^3W\tens \A)^G$ est 
aisé. Il suffit d'appliquer plusieurs fois la proposition 
\ref{propinv}. On tient compte de $\Stab\{1\}=\sigma_1\times\sigma_{m-1}$, $\Stab\{(12)\}=id_{(12)}\times\sigma_{m-2}$ et $\Stab\{1,2,3\}=\sigma_3\times\sigma_{m-3}$ et  on obtient
\begin{eqnarray}
\label{decomp1}
\H^0(\es^3W\tens\A)^{\sigm}&=&\H^0(L_1^{\tens 3}\tens 
A_1\tens\cdots\tens A_m)^{\sigma_1\times\sigma_{m-1}} \nonumber\\
&&\oplus\H^0(L_1^{\tens 2}\tens A_1\tens L_2\tens A_2\tens 
A_3\tens\cdots\tens 
A_m)^{id_{(12)}\times\sigma_{m-2}}\nonumber\\
&&\oplus\H^0(L_1\tens A_1\tens L_2\tens A_2\tens L_3\tens A_3\tens\cdots\tens 
A_m)^{\sigma_3\times\sigma_{m-3}}\nonumber\\
&=&\H^0(L^{\tens 3}\tens A)\tens\es^{m-1}(\H^0(A))\nonumber\\
&&\oplus\H^0(L^{\tens 2}\tens A)\tens\H^0(L\tens A)\tens\es^{m-2}(\H^0(A))\nonumber\\
&&\oplus\es^3\H^0(L\tens A)\tens\es^{m-3}(\H^0(A))
\end{eqnarray}

En restriction à l'ouvert $\U_{12}$ on a un \iso:
$$J^1L^3\tens A^{\tens 2}\boxtimes(A_3\tens\cdots\tens A_m)\oplus 
L_{\Delta}^2\tens A^{\tens 2}\boxtimes(\oplus_{i\ge 
3}L_i\tens(A_3\tens\cdots\tens A_m))\stackrel{\simeq}{\to}\gr_1\es^3W\tens \A$$
qui est dû au fait que $\gr_1(\es^3W\tens\A)=\gr_1(\es^3W)\tens\A$.
En utilisant la même proposition \ref{propinv} et en tenant compte de $\Stab\{3\}=\sigma_{m-3}$ on obtient 
\begin{eqnarray}
\label{decomp2}
\H^0(\gr_1\es^3W\tens\A)^{\sigm}&=&\H^0(\gr_1\es^3W\tens\A|_{\U_{12}})^{\sigma_2
\times\sigma_{m-2}}=\nonumber \\
&=&\left[\H^0(J^1L^3\tens A^{\tens
    2})\tens\es^{m-2}(\H^0(A))\right]\nonumber\\
&&\oplus\left[\H^0(L^{\tens 2}\tens A^{\tens 2})\tens\H^0(L_3\tens A_3\tens\cdots\tens 
A_m)^{\sigma_{m-3}}\right]\nonumber \\
&=&\left[\H^0(J^1L^3\tens A^{\tens
    2})\tens\es^{m-2}(\H^0(A))\right]\nonumber \\
&&\oplus\left[\H^0(L^{\tens 2}\tens A^{\tens 2})\tens\H^0(L\tens A)\tens\es^{m-3}(\H^0(A))\right]
\end{eqnarray}

\begin{rem}
{\rm Une section rationnelle est une section régulière sur un ouvert partout dense. Donc  $\nabla$ et $D$ se prolongent de manière évidente  aux
sections rationnelles du fibré $L^{\tens k}$ et la formule
$\nabla(sf)=df\tens s+f\nabla(s)$  est vraie pour $f$ fonction rationnelle sur
$X$, et $s$ section rationnelle de $L^{\tens k}$.}
\end{rem}

\begin{prop}

\label{matrice}

La matrice du \mor canonique
$$\alpha:\H^0(\es^3W\tens\A)^{\sigm}\to\H^0(\gr_1\es^3W\tens\A)^{\sigm}$$
dans les décompositions ci-dessus est de la forme
$$\left(
\begin{array}{ccc}
\widetilde{\nabla}&\widetilde{D}&0\\
0&\rho&\widetilde{\nu}\\
\end{array}
\right)$$

Les \mors $\rho$ et $\widetilde{\nu}$ sont 
$\O(X)$-linéaires et 
caractérisés par
$$\rho(s\tens t\tens a^{\tens_{\comp}(m-2)})=2sa\tens t\tens 
a^{\tens_{\comp}(m-3)}$$
lorsque $s\in\H^0(L^{\tens 2}\tens A)$, $t\in\H^0(L\tens A)$ et $a\in\H^0(A)$,
$$\widetilde{\nu}=-2\nu\tens id_{\es^{m-3}(\H^0(A))}.$$
où $\nu$ est l'opérateur défini au \ref{defnu} relatif à $L\tens A$.

Enfin le \mor $\widetilde{D}$ est caractérisé, pour $s$ section rationnelle de $L^{\tens 2}$, 
$t$ section rationnelle de $L$ et $a\in\H^0(A)$ tel que $sa\in\H^0(L^{\tens 2}\tens A)$ et $ta\in\H^0(L\tens A)$, par la formule
$$\tD(sa\tens_{\comp}ta\tens_{\comp}a^{\tens_{\comp}(m-2)})=D(s,t)a^2\tens_{\comp}a^{
\tens_{\comp}(m-2)}$$
et le \mor $\tna:\H^0(L^{\tens 
3}\tens A)\tens_{\comp}\es^{m-1}\H^0(A)\to\H^0(\gr_1(\es^3W_{12})\tens A^{\tens 
2})\tens_{\comp}\es^{m-2}\H^0(A)$ est caractérisé par:
$$\tna(sa\tens_{\comp}a^{\tens_{\comp}(m-1)})=-\nabla(s)a^2\tens 
a^{\tens_{\comp}(m-2)}$$
où $s$ est une section rationnelle de $L^{\tens 3}$  et 
$a\in\H^0(A)$ tel que $sa\in\H^0(L^{\tens 3}\tens A)$.

\end{prop}
\begin{rem}
\label{rem:les}
{\rm Les sections particulières qu'on a considérées pour décrire $\tD$ et $\tna$ sont  des générateurs des espaces vectoriels sur lesquels sont définis  ces morphismes.
Comme les \mors $\tD$ et $\tna$ existent ils sont caractérisés par l'image de ces générateurs. Puisqu'elles proviennent de \mors bien définis, les images de ces générateurs, qui sont a priori des sections rationnelles, sont bien des sections régulières.

À partir des expressions données, on peut facilement déduire les
expressions de $\tD$ et $\tna$ sur des sections de la forme
\hbox{$S\tens_{\comp}T\tens_{\comp}
  a^{\tens_{\comp}(m-2)}=\frac{S}{a}a\tens_{\comp}\frac{T}{a}a\tens_{\comp} a^{\tens_{\comp}(m-2)}$} ou \hbox{$S\tens_{\comp}a^{\tens_{\comp}(m-1)}=\frac{S}{a}a\tens_{\comp}a^{\tens_{\comp}(m-1)}$} en utilisant les propriétés des opérateurs $\nabla$ et $D$, pour obtenir ensuite les expressions sur des sections générales par polarisation, comme dans la remarque \ref{rem:par}. Ces expressions seront utilisées dans les lemmes \ref{morbeta} et \ref{morbeta2}.}
\end{rem}

\par{\bf Preuve de la proposition:}

Pour la première colonne soit $\gamma\in\H^0(L^{\tens 3}\tens A)\tens 
\es^{m-1}\H^0(A).$ D'après la remarque il suffit de traiter le cas où $\gamma=sa\tens_{\comp}a^{\tens_{\comp}(m-1)}$ avec pour $s$  une section 
rationnelle de $L^{\tens 3}$  et $a\in\H^0(A)$ tel que
$sa\in\H^0(L^{\tens 3}\tens A)$. Il suffit de prouver l'égalité de
l'énoncé sur l'ouvert $U\times U\times\cdots\times U$, $U$ étant un
ouvert de $X$ où $s$ est régulière et $a$ est inversible.
 La section invariante de $\es^3W\tens\A $ correspondante est alors

$S=s_1a_1\tens a_2\cdots a_m+s_2a_2\tens 
a_1a_3\cdots a_m+\cdots+s_ma_m\tens a_1\cdots a_{m-1}$
où $s_i=pr_i^*(s)$ et $a_i=pr_i^*(a)$. Comme les \fibs 
considérés sont inversibles les produits tensoriels s'identifient aux
produits symétriques et on peut écrire cette section comme 
$(s_1+s_2+\cdots+s_m)a_1a_2\cdots a_m$. L'image de $S$ dans 
$\F^1\es^3W\tens\A|_{\U_{12}} $ sera alors $(s_1+s_2)a_1\cdots 
a_m+(s_3+\cdots+s_m)a_1\cdots a_m$. Le deuxième terme appartient à 
$\F^2\es^3W\tens\A|_{\U_{12}} $ et l'image du premier terme modulo 
$\F^2\es^3W_{12}\tens\A $ est $-\nabla(s)a_1a_2\tens a_3\cdots a_m$ dans 
$\gr_1(\es^3W_{12})\tens A_1\tens A_2\boxtimes A_3\tens\cdots\tens A_m $ 
soit $-\nabla(s)a^2\tens a_3\cdots a_m$ dans $\gr_1(\es^3W_{12})\tens 
A_{\Delta}^2\boxtimes A_3\tens\cdots\tens A_m $ (car 
$A_1|_{\Delta_{12}}= A_2|_{\Delta_{12}}=A_{\Delta}$) ou encore 
$-\nabla(s)a^2\tens_{\comp}a^{\tens_{\comp}(m-2)}$ dans 
$\H^0(\gr_1(\es^3W_{12})\tens A^{\tens 2})\tens_{\comp}\es^{(m-2)}\H^0(A)$. D'où la première colonne de la matrice comme $(\tna,0)$.

\vs On procède de la même manière pour les autres colonnes. Soient
$a\in\H^0(A)$,  $s$ une section 
rationnelle de $L^{\tens 2}$ et $t$
une section 
rationnelle de $L$, telles que $sa\in\H^0(L^{\tens 2}\tens A)$ et
$ta\in\H^0(L\tens A)$. Partons d'une 
section de la forme $sa\tens ta\tens a^{\tens_{\comp}(m-2)}$ de 
$\H^0(L^{\tens 2}\tens A)\tens_{\comp}\H^0(L\tens 
A)\tens_{\comp}\es^{m-2}\H^0(A)$ comme dans l'énoncé. On se restreint
à un ouvert $U\times U\times\cdots\times U$, $U$ étant un ouvert de $X$ où $s$ et $t$ sont régulières et $a$ est inversible. La section 
invariante correspondante s'écrit
\begin{eqnarray*}
\left[(s_1t_2+s_2t_1)\prod_i a_i\right]+\left[(s_1+s_2)\sum_{i=3}^m t_i\prod_i 
a_i\right]+\left[(t_1+t_2)\sum_{i=3}^m s_i\prod_i a_i\right]+\left[\sum_{3\le i<j\le m}(s_it_j+s_jt_i)\prod_i a_i\right].\\
\end{eqnarray*}

\vs Le premier terme appartient à $\F^1\es^3W_{12}\tens\A $, le deuxième 
à $\es^2W_{12}\tens W^{12}\tens\A $, le troisième à 
$W_{12}\tens\es^2W^{12}\tens\A $, et le dernier à $\es^3W^{12}\tens\A $, 
donc son image dans 
$$\gr_1\es^3W_{12}\tens A^2\boxtimes 
A_3\tens\cdots\tens A_m \oplus L_{\Delta}^{\tens 2}\boxtimes L_3\tens 
A_3\tens\cdots\tens A_m $$ 
est
$$D(s,t)a^2\tens a_3\cdots a_m+2sa^2\tens t_3a_3\cdots a_m$$ 
soit
$$D(s,t)a^2\tens a^{m-2}+2Sa\tens T\tens a^{m-3}$$
dans
$$[\H^0(J^1L^3\tens A^{\tens
  2})\tens\es^{m-2}(\H^0(A))]\oplus[\H^0(L^{\tens 2}\tens A^{\tens
  2})\tens\H^0(L\tens A)\tens \es^{m-3}(\H^0(A))]$$ 
où $S=sa$, $T=ta$.

Soient $a\in\H^0(A)$,  $s, t, u$ des sections 
rationnelles de $L$  telles que $sa, ta, ua\in\H^0(L\tens A)$. On se
restreint à un ouvert $U\times U\times\cdots\times U$, $U$ étant un
ouvert de $X$ où $s$, $t$ et $u$ sont régulières et $a$ est inversible.
Si on part d'une section $sa\cdot ta\cdot ua\tens a^{m-3}$ de $\es^3\H^0(L\tens A)\tens\es^{m-3}\H^0(A)$, on obtient la section invariante~:
$$\left[(s_1t_2+s_2t_1)\sum_{i=3}^m u_i\prod_i a_i\right]+\left[(s_1u_2+s_2u_1)\sum_{i=3}^m 
t_i\prod_i a_i\right]+\left[(t_1u_2+t_2u_1)\sum_{i=3}^m s_i\prod_i a_i\right]+$$
$$+\left[(s_1+s_2)\sum_{3\le i\ne j\le m}t_iu_j\prod_i a_i\right]+\left[(t_1+t_2)\sum_{3\le 
i\ne j\le m}s_iu_j\prod_i a_i\right] +\left[(u_1+u_2)\sum_{3\le i\ne j\le 
m}t_is_j\prod_i a_i\right] +$$
$$+\left[\sum_{3\le i\ne j\le m,i\ne k}s_it_ju_k\prod_i a_i\right].$$

La classe dans $\gr_1(\es^3W\tens\A)|_{\U_{12}} $ de sa restriction à 
$\U_{12}$ est
$$[(s_1t_2+s_2t_1)a_1a_2\tens(u_3+\cdots+u_m)a_3\cdots 
a_m]+[(s_1u_2+s_2u_1)a_1a_2\tens(t_3+\cdots+t_m)a_3\cdots 
a_m]$$
$$+[(u_1t_2+u_2t_1)a_1a_2\tens(s_3+\cdots+s_m)a_3\cdots a_m]$$
et chacune de ces composantes appartient à $\gr_0(\es^2W_{12}\tens 
A)\tens W^{12}\tens A_3\tens\cdots\tens A_m $. Le reste est dans 
$\F^2\es^3W\tens\A|_{\U_{12}} $.

Comme on l'a déjà vu, l'image dans $\H^0(\U,L^{\tens 2}\tens 
A^{\tens 2})\tens\H^0(\U,L\tens A)\tens \es^{m-3}\H^0(\U,A)$ est 
$-2sata\tens ua\tens a^{m-3}-2saua\tens ta\tens a^{m-3}-2uata\tens 
sa\tens a^{m-3}$ donc si $S=sa$, $T=ta$, $U=ua$,  $\tnu$ associe 
à $STU\tens a^{m-3}\mapsto -2(ST\tens U\tens a^{m-3}+SU\tens T\tens 
a^{m-3}+UT\tens S\tens a^{m-3})$, et la troisième colonne de la 
matrice s'écrit $(0,\tnu)$.$\Box$

\begin{rem}
\label{rem:par}
{\rm Par polarisation on peut trouver l'expression de $\rho$ sur des sections 
différentes, où $\check{a_i}$ signifie qu'on omet le terme $a_i$ de l'expression~:

$\rho(S\tens T\tens a_3a_4\cdots a_m)=\frac{1}{m-2}\sum_{i=3}^m 
Sa_i\tens T\tens a_3\cdots \check{a_i}\cdots a_m.$}\end{rem}


\subsection{Sections de $\es^3\maV_3\tens\D$}
\label{5.6}

On prend $L=\O(3)$, et $A=\O(1)$ sur $X=\pp$. Alors
\begin{eqnarray*}
L^{\tens 3}\tens A&=&\O(10)\\
L^{\tens 3}\tens A^{\tens 2}&=&\O(11)\\
L^{\tens 2}\tens A&=&\O(7)\\
L\tens A&=&\O(4)\\
L^{\tens 2}\tens A^{\tens 2}&=&\O(8).
\end{eqnarray*}

On pose $E=\H^2(\pp,\O(1))$ et on suppose que $m\ge 13$ (afin que 
$a=m-2\ge b=11$). Pour nous $m=n+9$ avec $n\ge 11$. Cela marche aussi pour $n\ge 6$, et cela donne deux fa\c cons de faire le calcul pour $6\le n\le 11$. 

On veut calculer le noyau du \mor $\alpha$:
$$\left[\H^0(L^{\tens 3}\tens A)\tens\es^{m-1}(\H^0(A))\right]\oplus\left[\H^0(L^{\tens 2}\tens A)\tens\H^0(L\tens A)\tens\es^{m-2}(\H^0(A))\right]
\oplus\left[\es^3\H^0(L\tens A)\tens\es^{m-3}(\H^0(A))\right]$$
$$\to\left[\H^0(J^1L^3\tens A^{\tens 2})\tens\es^{m-2}(\H^0(A))\right]\oplus\left[\H^0(L^{\tens 2}\tens A^{\tens 2})\tens\H^0(L\tens A)\tens\es^{m-3}(\H^0(A))\right]$$
donné par la matrice de la proposition \ref{matrice}:
$$\left(
\begin{array}{ccc}
\tna&\tD&0\\
0&\rho&\tnu\\
\end{array}
\right)
$$

Ici $\H^0(\Omega^1L^3\tens A^{\tens 2})=\es^{10,1}E$ se calcule à 
partir de la suite exacte d'Euler. Le \mor $\H^0(\gr_1(\es^3W_{12})\tens A^{\tens 
2})\to\H^0(L^{\tens 3}\tens A^{\tens 2})=\es^{11}E$ est un \mor 
surjectif puisque non nul et que $\H^0(L^{\tens 3}\tens A^{\tens 2})$ est une 
\rep irréductible, d'où un scindage de la suite exacte
$$0\to\H^0(\Omega^1L^3\tens A^{\tens 2})\to\H^0(\gr_1(\es^3W_{12})\tens 
A^{\tens 2})\to\H^0(L^{\tens 3}\tens A^{\tens 2})\to 0.$$
Cela entraîne $\H^0(\gr_1(\es^3W_{12})\tens A^{\tens
  2})=\es^{10}E\tens E$.

Avant le résultat final on va donner deux lemmes préliminaires~:

\begin{lemme}

Le \mor $D'$
$$\begin{array}{ccccccc}
\H^0(L^{\tens 2}\tens A)\tens_{\comp}\H^0(L\tens A) & = &\es^{7}E\tens 
\es^4E  & \to &  \H^0(\gr_1(\es^3W_{12})\tens A^{\tens 2}) & = &\es^{10}E\tens E\\
&&sa\tens ta&\mapsto&D(s,t)a^2&&
\end{array}$$
est conjugué au \mor composé des contractions naturelles.

\end{lemme}

\par{\bf Preuve du lemme:}

Il suffit de montrer que  $D'$ est surjectif. On a vu 
que $\epsilon_+(D(s,t))=-2st$ donc $\epsilon_+(D(s,t)a^2)=-2sta^2$, 
donc le \mor composé $\epsilon_+\circ D'$ de $\H^0(L^{\tens 2}\tens 
A)\tens_{\comp}\H^0(L\tens A)=\es^{7}E\tens\es^4 E$ dans $\H^0(L^{\tens 
3}\tens A^{\tens 2})=\es^{11}E$ est conjugué au \mor de multiplication 
des sections, et il est par suite non nul.

Si on se place sur l'ouvert $\U$ défini par $X\ne 0$, $Y\ne 0$ et on 
considère les sections $s=Z^6$, $t=Z^3$, $a=X$, $fa=Y$ avec 
$f=\frac{Y}{X}$, comme 
$$D'(sa\tens tfa)=D(s,tf)a^2=fD(s,t)a^2-st{\rm d}fa^2$$ 
et
$$D'(sfa\tens ta)=D(fs,t)a^2=fD(s,t)a^2+st{\rm d}fa^2,$$
on obtient que
$$D'(sfa\tens ta-sa\tens tfa)=2st{\rm d}fa^2.$$
Ces sections sont en effet globales et ceci montre que
$D'(Z^6Y\tens Z^3X-Z^6X\tens Z^3Y)$ est non nul et appartient à 
$\H^0(\Omega^1L^3\tens A^{\tens 2})$, vu comme sous-espace de 
$\H^0(\gr_1(\es^3W)\tens A^{\tens 2})$.
Le \mor $D'$ est donc non nul sur chacune des composantes \irrs de la \rep 
$\es^{10}E\tens E$, donc il est surjectif.$\Box$

\begin{lemme}
\label{rho}
L'image de $\Ker\,\tD$ par le \mor $\rho$ est incluse dans l'image de $\tnu$.
\end{lemme}

\par{\bf Preuve du lemme:}
On rappelle que $\tD=D'\tens id$.
On sait que $D'$ est conjugué au \mor composé des contractions 
naturelles $\es^7E\tens \es^4E\to\es^{10}E\tens E$. On a un diagramme 
commutatif:
$$\begin{array}{ccc}
{\es^7E\tens \es^4E\tens\es^{m-2}E} & \rightto{\rho} & {\es^8E\tens \es^4E\tens\es^{
m-3}E}\\
\downto{D'\tens id}& &\downto{D''\tens id}\\
{\es^{10}E\tens E\tens\es^{m-2}E}&\rightto{\rho'}&{\es^{11}E\tens 
E\tens\es^{m-3}E}
\end{array}$$
dans lequel les flèches verticales $D'$ et $D''$ sont les composées 
de contractions naturelles entre premier et second facteur, et les 
flèches horizontales $\rho$ et $\rho'$ sont des contractions entre 
premier et troisième facteur.

Il en résulte que $\rho(\Ker\, D'\tens id)$ est contenu dans 
$\Ker(D''\tens id)$. Ce dernier noyau est évidemment 
$(\es^{10,2}E+\es^{9,3}E+\es^{8,4}E)\tens\es^{m-3}E$ lequel est bien contenu 
dans l'image de $\tnu$, d'après le lemme \ref{tnu} de la section préliminaire.$\Box$

\begin{prop}

\label{propo}

\begin{itemize}

\item {\rm (i)} L'espace des sections 
$\Ker\,\alpha=\H^0(\es^3\maV_3\tens\de)$ est isomorphe à 
$\Ker(\tna,\tD)\oplus\Ker\,\tnu$.

\item {\rm (ii)} Sur l'ouvert $\hils$ considéré, 
$\coker\,\alpha=\H^1(\es^3\maV_3\tens\de)$ est isomorphe à 
$\coker(\tna,\tD)\oplus\coker\,\tnu$.

\end{itemize}

\end{prop}

\par{\bf Preuve:}

Cela revient à vérifier que dans la suite exacte du serpent 
associé au diagramme
$$\begin{array}{ccccccccc}
{0}&\rightto{}&{\mathfrak{C}}&\rightto{}&{\mathfrak{A}\oplus\mathfrak{B}\oplus \mathfrak{C}}&\rightto{}&{\mathfrak{A}\oplus 
\mathfrak{B}}&\rightto{}&{0}\\
&&\downto{\tnu}&&\downto{\alpha}&&\downto{(\tna,\tD)}&&\cr
{0}&\rightto{}&{\mathfrak{E}}&\rightto{}&{\mathfrak{D}\oplus\mathfrak{E}}&\rightto{}&{\mathfrak{D}}&\rightto{}&{0}
\end{array}$$
avec ${\mathfrak{A}}=\es^{10}E\tens\es^{m-1}E$, ${\mathfrak{B}}=\es^7E\tens\es^4E\tens\es^{m-2}E$, 
${\mathfrak{C}}=\es^3(\es^4E)\tens\es^{m-3}E$, ${\mathfrak{D}}=\es^{10}E\tens E\tens\es^{m-2}E$, 
${\mathfrak{E}}=\es^8E\tens\es^4E\tens\es^{m-3}E$, c'est à dire dans la suite:
$$0\to\Ker\,\tnu\to\Ker\,\alpha\to\Ker(\tna,\tD)\stackrel{(0,\rho)}{\to}\coker\,\tnu\to \coker\,\alpha\to\coker(\tna,\tD)\to 0$$
le \mor de liaison $(0,\rho)$ est nul. Cela revient à vérifier que 
l'image de $\Ker(\tna,\tD)$ par $(0,\rho)$ est contenue dans l'image de $\tnu$.

On dispose également d'une suite exacte:
$$\begin{array}{ccccccccc}
0&\to&\Ker\,\tD&\stackrel{\gamma}{\to}&\Ker(\tna,\tD)&\stackrel{\beta}{\to}&{\mathfrak{A}}& \\
&v&\mapsto&(0,v)&&&&&\\
&&&(u,v)&\mapsto&u&&&
\end{array}$$

\begin{lemme}
\label{morbeta}
Le \mor $\beta$ est surjectif.
\end{lemme}

\par{\bf Preuve du lemme:}
 
Pour cette preuve on peut en effet montrer seulement que les éléments $u$ de la forme $\omega^3a\tens a_2\cdots a_m$ ont un antécédent par $\beta$, où $\omega\in\H^0(L)$ et $a, a_2, \ldots,a_m\in\H^0(A)$ (on rappelle que $L=\O(3)$ et $A=\O(1)$ mais on préfère travailler ici avec les \fibs quelquonques $L$ et $A$ de départ pour une meilleure compréhension) et ceci puisque sur $\pp$ le \mor naturel $\es^3\H^0(L)\tens\H^0(A)\tens\es^{m-1}\H^0(A)\to\H^0(L^{\tens 3}\tens A)\tens\es^{m-1}\H^0(A)$ est surjectif et les éléments considérés engendrent $\es^3\H^0(L)\tens\H^0(A)\tens\es^{m-1}\H^0(A)$.
On prend $v=\sum_{i=2}^{m}(2\omega^2a\tens_{\comp}\omega
a_i+\omega^2a_i\tens_{\comp}\omega a)\tens a_2\cdots\check{a_i}\cdots
a_m$ et on remarque que $(u,-v)\in\Ker(\tna,\tD)$ est une préimage de
$u$ par $\beta$. Rappelons que la notation $\check{a_i}$ signifie
qu'on omet le terme $a_i$ de l'expression. On utilise les formules de $\tna$ et $\tD$ sur des
sections différentes déduites de celles données dans l'énoncé de la
prop. \ref{matrice}, par polarisation, comme dans les remarques
\ref{rem:les} et \ref{rem:par}:
$$\tna(f\tens a_2\cdots a_m)=\sum_{i=2}^m\nabla(\frac{f}{a_i})a_i^2\tens a_2\cdots\check{a_i}\cdots a_m$$
et
$$\tD(s\tens t\tens a_3\cdots a_m)=\nabla(\frac{s}{t})t^2\tens a_3\cdots a_m$$
et le calcul nous donne que $\tna(u)+\tD(-v)=\tna(u)-\tna(u)=0$.$\Box$

\begin{lemme}
\label{morbeta2}
L'image   par le \mor $(0,\rho)$ de la préimage par $\beta$ d'un
élément de ${\mathfrak{A}}$ est
contenue dans 
l'image de $\tnu$.
\end{lemme}

\par{\bf Preuve du lemme:}

Dans le lemme précédent on a trouvé un antécédent $v$ pour chaque
générateur de ${\mathfrak{A}}$. Comme deux antécédents d'un 
élément $u$ de $\mathfrak{A}$ diffèrent par un élément de $\Im\gamma$,
et que  $(0,\rho)(\Im\gamma)\subset\Im\tnu$  (lemme \ref{rho}),
il suffit  de trouver un antécédent $w$ par $\tnu$ de la section
$\rho(v)=\sum_{i\ne j}(2\omega^2aa_j\tens\omega
a_i+\omega^2a_ia_j\tens\omega a)\tens a_2\cdots
\check{a_i}\cdots\check{a_j}\cdots a_m$. On prend $w=\sum_{i\ne
  j}(\omega a\cdot\omega a_i\cdot\omega a_j)\tens a_2\cdots
\check{a_i}\cdots\check{a_j}\cdots a_m$ et on vérifie bien que
$\tnu(w)=\rho(v)$.$\Box$

Ceci montre que $(0,\rho)(\Ker(\tna,\tD))$.$\Box$

\subsection{Calcul final}

Le résultat de la proposition \ref{propo} nous donne 
\begin{eqnarray*}
\dim\H^0(\es^3\V_3\tens\de)&=&\dim{\mathfrak{A}}+\dim\Ker\tD+\dim\Ker\tnu=\\
&=& \dim\es^{10}E\tens\es^{n+8}E+\dim\es^7E\tens\es^4E\tens\es^{n+7}E-\dim\es^{10}E\tens E\tens\es^{n+7}E+\\
&+&\dim\es^{n+6}E(\dim\es^{6,6}E+\dim\es^{7,4,1}E+\dim\es^{8,2,2}E+\dim\es^{6,4,2}E+\dim\es^{4,4,4}E).
\end{eqnarray*}
 
À partir de ce résultat et des annulations de la \coh supérieure 
indiquées dans la démonstration du corollaire \ref{h1v} remplacées dans la suite spectrale associée 
à la résolution (\ref{**}) de $\es^3\R\tens\de$, on obtient une suite 
exacte
$$0\to\Lambda^3(\es^3E)\tens\H^0(\de)\to\Lambda^2(\es^3E)\tens\H^0(\V_3\tens\de) \to \es^3E\tens\H^0(\es^2\V_3\tens\de)\to $$
$$\to\H^0(\es^3\V_3\tens\de)\to\H^0(\es^3\R\tens\de)\to 0$$
D'où encore
\begin{eqnarray*}
\dim\H^0(\es^3\R\tens\de)&=&\dim\H^0(\es^3\V_3\tens\de)-
10(\dim\es^7E\tens\es^{n+8}E+\dim\es^2(\es^4E)\tens\es^{n+7}E)\\
&-&\dim\es^8E\tens\es^{n+7}E+45\dim\es^4E\tens\es^{n+8}E-120\dim\es^{n+9}E\\
&=&\frac{1}{2}(n+1)(n+2)
\end{eqnarray*}
et cela pour tout $n$ tel que $6\le n\le 19$. On voit que dans ce cas on couvre  aussi le résultat obtenu au corollaire \ref{h1v} pour $l=2$.


\section{Conclusion}

Pour étendre les résultats ci-dessus au cas $n\ge 20$, on a besoin

1) d'étendre le \th d'annulation de la cohomologie supérieure des \fibs $\es^k\V\tens\de$ sur le schéma de Hilbert $\hil$;

2) de faire intervenir $\H^0(\gr_i(\es^lW)\tens\de)$ (pour $i\ge 2$)
pour le calcul de $\H^0(\es^l\vkde)$, pour $l\ge 3$.

\par{\bf Remerciements~:} Les explications détaillées de J. Le Potier
et l'aide de N. Dan ont rendu possible ce texte. 
Une faute grave
dans le texte initial m'a été signalée par C. Mourougane. Je les
remercie, ainsi que D. Roessler, pour leur lecture patiente et les
nombreuses corrections.



\end{document}